\documentclass[preprint,11pt]{article}
\usepackage{amssymb}
\usepackage{mathrsfs}
\usepackage{amsfonts}
\usepackage{graphicx}
\usepackage{colortbl,dcolumn}
\usepackage{tabularx}
\usepackage{amsmath}
\usepackage{psfrag}
\usepackage{booktabs}
\usepackage{array}
\usepackage{cite}
\numberwithin{equation}{section}
\usepackage[top=1.2in, bottom=1.2in, left=0.9in, right=0.9in]{geometry}

\newcommand{\dd}{\text{d}}
\newtheorem{theorem}{Theorem}[section]
\newtheorem{remark}[theorem]{Remark}
\newtheorem{assumption}[theorem]{Assumption}

\newtheorem{lemma}[theorem]{Lemma}

\begin{document}

%
 \title {Error estimates of finite element method for semi-linear stochastic strongly damped wave equation
\footnote{R.Q. was supported by Research Fund for Northeastern University at Qinhuangdao (No.XNB201429),
Fundamental Research Funds for Central Universities  (No.N130323015),  Science and Technology Research
Funds for Colleges and Universities in Hebei Province (No.Z2014040),  Natural Science Foundation of Hebei
Province (No.A2015501102). X.W. was supported by NSF of China (No.11301550,  No.11571373) and Innovation Program of Central South University.
%
}
}

\author{
Ruisheng Qi$\,^\text{a}$,  \quad Xiaojie Wang$\,^\text{b}$ \\
\footnotesize $\,^\text{a}$ School of Mathematics and Statistics, Northeastern University at Qinhuangdao, Qinhuangdao, China\\
\footnotesize qirsh@neuq.edu.cn\; and \;qiruisheng123@sohu.com\\
\footnotesize $\,^\text{b}$ School of Mathematics and Statistics, Central South University, Changsha, China\\
\footnotesize x.j.wang7@csu.edu.cn\; and \;x.j.wang7@gmail.com
}
\maketitle
\begin{abstract}\hspace*{\fill}\\
  \normalsize
In this paper, we consider a semi-linear stochastic strongly damped wave equation driven by additive Gaussian noise. Following a semigroup framework, we establish existence, uniqueness and space-time regularity of
a mild solution to such equation. Unlike the usual stochastic wave equation without damping,
the underlying problem with space-time white noise ($Q = I$) allows for a mild solution with
a positive order of regularity in multiple spatial dimensions.
Further, we analyze  a spatio-temporal discretization of the problem,  performed by a standard finite element
method in space and a well-known linear implicit Euler scheme in time.
The analysis of the approximation error forces us to significantly enrich existing error estimates of semidiscrete
and fully discrete finite element methods  for the corresponding  linear deterministic equation.
The main results show optimal convergence rates in the sense that the orders of convergence in space
and in time coincide with the orders of the spatial and temporal regularity of the mild solution, respectively.
Numerical examples are finally included to confirm our theoretical findings.

  \textbf{\bf{Key words.}}
strongly damped wave equation, Wiener process, finite element method, linear implicit Euler scheme, strong approximation
\end{abstract}

\section{Introduction}

The present work is concerned with the following semi-linear stochastic evolution equation subject to additive noise, described by
\begin{align}\label{eq:intro-SSDVE}
  \left\{\begin{array}{ll}
\dd u_t = \alpha L u_t \, \dd t + L u \, \dd t  + F(u) \, \dd t +  \dd W(t), & \text{ in }\; \mathcal{D}\times (0,\; T],\\
 u(\cdot,0)=\varphi,\; u_t ( \cdot, 0 )=\psi,& \text{ in } \; \mathcal{D},\\
  u=0, & \text{ on } \; \partial \mathcal{D}\times (0,\; T],
  \end{array}\right.
\end{align}
where $\mathcal{D}\subset\mathbb{R}^d$, $d=1,2,3$, is a bounded, convex and polynomial domain with
a  boundary $\partial\mathcal{D}$ and $\alpha > 0$ is a fixed positive constant.
Let $L = \sum_{i, j = 1}^d \tfrac{\partial}{ \partial x_i}  \big(  l_{ij} (x) \tfrac{\partial}{ \partial x_j}  \big),
x \in \mathcal{D} $ be a linear second-order elliptic operator with smooth coefficients and $\{l_{ij}\}$
 being uniformly positive definite.  Let $\{W(t)\}_{t\in [0, T]} $ be a (possibly cylindrical)
 $Q$-Wiener process on $\big( L_2(\mathcal{D}), \| \cdot \|, (\cdot, \cdot) \big)$, defined on a stochastic basis
$(\Omega,\mathcal{F},\mathbb{P},\{\mathcal{F}_t\}_{t \in [0, T]})$ with
respect to a normal filtration $\{\mathcal{F}_t\}_{t\in [0, T] } $ and let
$\varphi, \psi$ be $\mathcal{F}_0$-measurable random variables.

The considered problem \eqref{eq:intro-SSDVE} is referred to as stochastic strongly damped wave equation (SSDWE for short) thereafter.
The deterministic counterpart of \eqref{eq:intro-SSDVE} finds many applications in viscoelastic theory \cite{fitzgibbon1981strongly,massatt1983limiting,pata2005strongly}, and its linear version has been numerically studied by \cite{larsson1991finite,thomee2004maximum}, where a finite element method is used for spatial discretization
and rational approximations for analytic semigroup. Particularly when $\alpha = 0$,  the problem \eqref{eq:intro-SSDVE} reduces to a stochastic wave equation (SWE) without damping,  numerical approximations of which have been recently studied by many authors \cite{anton2015full,cao2007spectral,cohen2013trigonometric,cohen2015fully,hausenblas2010weak,kovacs2012BITweak,kovacs2013BITweak,walsh2006numerical,quer2006space,wang2014higher,wang2015exponential,jacobe2015weak,jiang2015stochastic,qi2013weak,qi2013full}.
In contrast to the SWE case ($\alpha = 0$), the stochastic strongly damped wave equations ($\alpha > 0$)
are much less well-understood, from both theoretical and numerical point of view.
In \cite[Example 6.25]{da2014stochastic}, a linear version of SSDWE with multiplicative noise was
examined and its unique mild solution was verified. To the best of our knowledge, regularity analysis
and numerical treatment of such stochastic problem are both missing in the literature.
This article aims to fill the gap and investigate the regularity properties and strong approximations of SSDWE like \eqref{eq:intro-SSDVE}.
%

Reformulating \eqref{eq:intro-SSDVE} as a Cauchy problem of first order in a Hilbert space, we follow the
semigroup framework as in \cite{da2014stochastic} to show existence, uniqueness and space-time
regularity of a mild solution to \eqref{eq:intro-SSDVE}. Under some standard assumptions
(Assumptions \ref{ass:F}-\ref{ass:initial-value}),  it is revealed that (see Theorem \ref{thm:main-regularity}),
the unique mild solution $\{ u(t)\}_{t \in [0, T]}$ exhibits the Sobolev and H\"older regularity properties as follows,
\begin{equation} \label{eq:intro-u-regularity}
u \in  L^{\infty} \big( [0, T]; L^2( \Omega; \dot{H}^{\gamma+1} ) \big), \quad
\sup_{ t, s \in [0, T], t \neq s } \tfrac{ \| u(t) - u(s) \|_{L^2(\Omega;\dot{H}^0)} } { |t-s|^{(\gamma + 1)/2 } }
< \infty,
\end{equation}
where, as specified later, $\dot{H}^{\delta} := D(A^{\frac{\delta}{2} }) $, $\delta \subset\mathbb{R}$ and
the parameter $\gamma \in [-1, 1]$ satisfying
$\|A^{\frac{\gamma-1}{2}}Q^{\frac{1}{2}}\|_{\mathrm{HS}} < \infty$ quantifies the spatial correlation of
the noise process (Assumption \ref{ass:QWiener}). If $\gamma \in [0, 1]$ and $\psi  \in L^2( \Omega; \dot{H}^{\gamma} )$, then
\begin{equation}\label{eq:intro-v-regularity}
u_t \in  L^{\infty} \big( [0, T];  L^2( \Omega; \dot{H}^{\gamma} ) \big),
\quad
\sup_{ t, s \in [0, T], t\neq s } \tfrac{ \| u_t(t) - u_t(s) \|_{L^2(\Omega;\dot{H}^0)} } { |t-s|^{\gamma/2 } } < \infty.
\end{equation}
In order to achieve \eqref{eq:intro-u-regularity}-\eqref{eq:intro-v-regularity}, we exploit further spatial
and temporal regularity  properties of the linear deterministic equation (Lemmas
\ref{lem:linear-determin-spatial}-\ref{Lem:timez}),  based on an existing spatial regularity result (Lemma \ref{lem:S-GROUP-spatial-regularity})  in \cite{larsson1991finite}.  From \eqref{eq:intro-u-regularity}-\eqref{eq:intro-v-regularity}, it is easy to realize that, the mild solution of \eqref{eq:intro-SSDVE} ($\alpha > 0$) enjoys
higher spatial and temporal regularity than that of the usual stochastic wave equation ($\alpha = 0$),
which only admits a mild solution taking values in $L^{\infty} \big( [0, T];  L^2( \Omega; \dot{H}^{\gamma} \times \dot{H}^{\gamma-1} ) \big)$ and satisfying $\sup_{ t, s \in [0, T] } \tfrac{ \| u(t) - u(s) \|_{L^2(\Omega;\dot{H}^0)} } { |t-s|^ {\gamma} } < \infty$ under the same assumptions \cite{anton2015full,wang2015exponential}.
This benefits from smoothing effect  of the analytic semigroup $\mathcal{S}(t)$ generated by the dominant linear operator $\mathcal{A}$.
In particular, different from both the stochastic heat equation and the stochastic wave equation, the strongly damped problem driven by space-time white noise ($Q = I$) allows for a mild
solution with a positive order of regularity in multiple spatial dimensions ($d > 1$).
For example, the space-time white noise case when $d = 2$ admits a mild solution
$u \in L^{\infty} \big( [0, T]; L^2( \Omega, \dot{H}^{\alpha} ) \big)$ for any $\alpha < 1$
(consult Remark \ref{rem:space-time-white-noise} for more details).

As the second contribution of this article, we analyze the mean-square approximation errors
caused by finite element spatial semi-discretization and space-time full-discretization of \eqref{eq:intro-SSDVE}.
More precisely, we measure the discrepancy between the mild solution $(u(t), u_t(t))'$ and
the finite element spatial approximation $(u_h(t), u_{h,t}(t))'$ as follows (Theorem \ref{thm:conteh1}):
%
\begin{equation} \label{eq:intro-FEM-main}
\|u(t)-u_h(t)\|_{L^2(\Omega;\dot{H}^0)} = O \big( h^{1+\gamma} \big),
\quad
\|u_t(t)-u_{h,t}(t)\|_{L^2(\Omega;\dot{H}^0)} = O \big( h^{\gamma} \big),
\end{equation}
where the parameter $\gamma$ restricted to $\gamma \in [0, 1]$, similarly as before,  characterizes the spatial correlation of the Wiener process.
By a combination of the finite element method (FEM) together with a linear implicit Euler-Maruyama time-stepping scheme, we also investigate a spatio-temporal discretization of \eqref{eq:intro-SSDVE}.
As stated in Theorem \ref{thm:conteh2}, the corresponding strong approximation error satisfies
\begin{equation}
\label{eq:intro-fullFEM-main}
\|u(t_n)-U^n\|_{L^2(\Omega;\dot{H}^0)} = O \big( h^{\gamma+1}+k^{\frac{ \gamma+1 } {2} } \big),
\quad
\|u_t(t_n)-V^n\|_{L^2(\Omega;\dot{H}^0)} = O \big( h^{\gamma}+k^{\frac{\gamma}{2}} \big),
\: \gamma \in [0, 1].
\end{equation}
Here $U^n$ and $V^n$ are, respectively, full-discrete approximations of $u(t_n)$ and $u_t(t_n)$.
Comparing the convergence results \eqref{eq:intro-FEM-main}-\eqref{eq:intro-fullFEM-main}
with the regularity results \eqref{eq:intro-u-regularity}-\eqref{eq:intro-v-regularity},
one can readily observe that, the convergence rates obtained here are optimal in the sense that
the rates of convergence in space and in time coincide with the orders of the spatial and
temporal regularity of the mild solution, respectively. This essentially differs from the SWE setting ($\alpha = 0$),
where the strong rates $O(h^{ \frac23\gamma} + k^{\frac{\gamma} {2} })$ of
the FEM coupled with the linear implicit Euler scheme are lower than orders of the spatial and temporal regularity of the mild solution (e.g., \cite{kovacs2010finite,kovacs2013BITweak}).


Before proving \eqref{eq:intro-FEM-main}-\eqref{eq:intro-fullFEM-main},
we formulate in section \ref{sec:Error-estimate-determin}  a rich variety of error estimates for
the finite element semi-discretization and full-discretization of the corresponding
deterministic linear problem.
Some of such error estimates can be straightforwardly derived from existing ones in \cite{larsson1991finite}
by ingenious modifications or by interpolation arguments  (Theorems \ref{thm:error-FEM1},
\ref{thm:error-FEM3}). Nevertheless, we must stress that, error estimates available in \cite{larsson1991finite}
are far from enough for the purpose of our error analysis. For instance, as one can see later,
two completely new error estimates of integral form such as  \eqref{eq:error-FEM-u-v}
and \eqref{eq:error-full-u-v}  are indispensable in the error analysis and their proofs turn out to be quite involved.
%
To show the error estimate \eqref{eq:error-FEM-u-v} for the semi-discretization, we rely on energy arguments and interpolation theory (see the proof of Theorem \ref{thm:error-FEM2}).
The proof of \eqref{eq:error-full-u-v} for the full-discretization is, however, more complicated and more technical.
In addition to energy arguments and interpolation theory, we need some further integral versions of regularity results of the linear deterministic problem as presented in Lemma \ref{lem:appen1}.
%
Armed with these error estimates, we are then able to establish \eqref{eq:intro-FEM-main}-\eqref{eq:intro-fullFEM-main} for the stochastic problem (see section \ref{sec:stochastic-main-result} for the details).


The outline of this paper is as follows. In the next section, some preliminaries are collected and the well-posedness of the considered problem is elaborated. Section \ref{sec:Error-estimate-determin} is devoted to error estimates of semi-discrete and full-discrete finite element method for the corresponding deterministic linear problem. The main convergence results for the stochastic problem are presented in section \ref{sec:stochastic-main-result}. Numerical experiments are finally performed in section \ref{sec:numerical-exam} to confirm the theoretical results.

\section{The stochastic strongly damped wave equation}
\label{sec:regularity-stochastic-problem}
 Let $U$ and $H$
 be two separable $\mathbb{R}$-Hilbert spaces and by $\mathcal{L}(U, H)$ we denote the Banach space of all linear bounded operators from $U$ into $H$ and by $\mathcal{L}_2(U, H)$ the Hilbert space of all Hilbert-Schmidt operators from $U$ into $H$.
%
When $H = U$, we write $\mathcal{L}(U) =\mathcal{ L}(U,U)$
and $\text{HS} = \mathcal{L}_2(U,U)$ for ease of notation.
Also, we denote the space of the Hilbert-Schmidt operators from
$Q^{\frac{1}{2}}(U)$ to $H$ by
$\mathcal{L}_2^0:= \mathcal{L}_2 (Q^{\frac{1}{2}}(U),H)$ and the corresponding
norm is given by
$
\|\Gamma \|_{\mathcal{L}_2^0}
=\|\Gamma Q^{\frac{1}{2}}\|_{\mathcal{L}_2(U, H) }.
$
It is well-known that
$
\|ST\|_{\mathcal{L}_2(U,H)}\leq \|T\|_{\mathcal{L}_2(U,H)}\|S\|_{\mathcal{L}(U)},
T\in \mathcal{L}_2(U,H),\; S\in \mathcal{L}(U).
$
%
Additionally, we denote $A = - L$ with the domain $D(A) = H^2 (\mathcal{D}) \cap H_0^1 (\mathcal{D})$
and define
\[
\dot{H}^s=D(A^{\frac{s}{2}}),\;\|v\|_{\dot{H}^s}=\|A^{\frac{s}{2}} v\|=|v|_s,\;s\in \mathbb{R}.
\]
It is clear that $\dot{H}^0=L_2(\mathcal{D})$, $\dot{H}^1=H_0^1 (\mathcal{D}) $ and
$\dot{H}^2 = H^2 (\mathcal{D}) \cap H_0^1 (\mathcal{D})$.
Now we reformulate the stochastic equation (\ref{eq:intro-SSDVE}) as the following abstract form
\begin{align}\label{stojab}
\left\{\begin{array}{ll}
\dd X(t)=- \mathcal{A}X(t) \, \dd t + \mathbf{F} (X(t)) \, \dd t + \mathbf{B} \, \dd W(t), &
t \in (0, T],\\
X(0)=X_0,&
\end{array}\right.
\end{align}
where $X(t)=(u(t),u_t(t))'$, $X_0=(\varphi,\psi)'$ and
\[\mathcal{A}:=\biggl[\begin{array}{c
c}0&-I\\A&\alpha A\end{array}\biggr],\quad \mathbf{F}(X):=
\biggl[
      \begin{array}{c}0 \\ F(u)\end{array}
\biggr]\;\: \text{ and } \;
\mathbf{B}:=\biggl[\begin{array}{c}0\\I\end{array}\biggr].
\]
Here $-\mathcal{A}$ generates an analytic semigroup $\mathcal{S}(t) = e^{ - t \mathcal{A} } $ in
$\dot{H}^s \times \dot{H}^{ s-\sigma}$, $s \in \mathbb{R},  \sigma \in [0, 2]$
(see \cite[Lemma 2.1]{larsson1991finite}).
%
For the purpose of the existence, uniqueness and regularity of the mild solution to (\ref{stojab}),
we put standard assumptions on the nonlinear term $F$,  the noise process $W(t)$
and the initial data $(\varphi, \psi)'$.

\begin{assumption}\label{ass:F}
(Nonlinearity)
Let $F \colon \dot{H}^0  \rightarrow \dot{H}^0 $ be a deterministic mapping such that
\begin{align}
\|F(x)-F(y)\|&\leq K \|x-y\|,\quad\;\forall\; x, y \in \dot{H}^0,  \label{asssmad}\\
\|F(x)\|&\leq K (1 + \|x\| ), \quad\;  \forall\;  x \in \dot{H}^0,  \label{assumdd}
\end{align}
where $K\in (0, \infty)$ is a positive constant.
\end{assumption}
\begin{assumption}\label{ass:QWiener}
(Q-Wiener process)
Let $W(t)$ be a (possibly cylindrical) $Q$-Wiener process on $ \dot{H}^0 $,
with the covariance operator $Q \colon \dot{H}^0  \rightarrow  \dot{H}^0$ being a symmetric nonnegative operator  satisfying
\begin{align}\label{jiezz1}
\|A^{\frac{\gamma-1}{2}}Q^{\frac{1}{2}}\|_{\mathrm{HS}} < \infty, \quad  \text{ for some } \ \gamma \in [-1, 1].
\end{align}
\end{assumption}

\begin{assumption}\label{ass:initial-value}
(Initial data)
Let $\varphi , \psi $ be $\mathcal{F}_0$-measurable and
$(\varphi,\psi)' \in L^2 (\Omega; \dot{H}^{\gamma+1})\times  L^2(\Omega; \dot{H}^{\gamma-1})$.
\end{assumption}

Here we let $\textbf{E}$ be the expectation in the probability space and let $L^2(\Omega;H)$ be
the space of $H$-valued integrable random variables, equipped with the norm
$\|v\|_{L^2(\Omega;H)}= \big( \textbf{E} [\|v\|_H^2] \big)^{\frac{1}{2}}$.
%
Owing to the above assumptions, we have the following regularity results of the mild solution of (\ref{stojab}).
\begin{theorem}\label{thm:main-regularity}
Under Assumptions \ref{ass:F}-\ref{ass:initial-value}, the problem \eqref{stojab} admits a unique mild solution given by
\begin{align}\label{eq:mild-SPDE}
X(t)=\mathcal{S}(t)X_0+\int_0^t\mathcal{S}(t-s)\mathbf{F}(X(s)) \,\mathrm{d} s
+\int_0^t\mathcal{S}(t-s)\mathbf{B} \, \mathrm{d}W(s), \quad t \in [0, T], \quad a.s..
\end{align}
Furthermore, the mild solution $\{ u(t) \}_{t \in [0, T]}$ has the following space-time regularity properties
\begin{align}
\sup_{t \in [0, T]} \| u(t) \|_{ L^2( \Omega;\dot{H}^{\gamma+1} ) }
      & \leq C \big( 1 + \|\varphi\|_{L^2 (\Omega; \dot{H}^{\gamma+1}) }
          + \| \psi \|_{L^2(\Omega; \dot{H}^{\gamma - 1} )}
          \big),
\label{eq:thm1-u-spatial2} \\
\| u(t) - u(s) \|_{L^2(\Omega;\dot{H}^0)}
      & \leq C |t-s|^{\frac{\gamma + 1} {2} }
          \big(
          1 + \|\varphi\|_{ L^2(\Omega; \dot{H}^{\gamma+1} ) }
          +\|\psi\|_{ L^2(\Omega; \dot{H}^{\gamma -1} )}
           \big).
 \label{eq:thm1-u-temporal2}
\end{align}
Additionally, if Assumption \ref{ass:QWiener} is satisfied with  $\gamma \in [0, 1]$ and
$\psi  \in L^2( \Omega; \dot{H}^{\gamma} )$, then
\begin{align}
\sup_{t \in [0, T]} \| u_t(t) \|_{L^2(\Omega;\dot{H}^{\gamma})}
      & \leq C \big(
      1 + \|\varphi\|_{L^2(\Omega; \dot{H}^{ \gamma} )}
           + \|\psi\|_{L^2(\Omega; \dot{H}^{\gamma})}
          \big),
\label{eq:thm1-v-spatial2} \\
\| u_t (t) - u_t(s) \|_{L^2(\Omega;\dot{H}^0)}
      & \leq C |t-s|^{\frac{\gamma } {2} }
          \big(
          1 + \|\varphi\|_{ L^2(\Omega; \dot{H}^{\gamma} ) }
          +\|\psi\|_{ L^2(\Omega; \dot{H}^{\gamma} )}
          \big).
 \label{eq:thm1-v-temporal2}
\end{align}
\end{theorem}
\begin{remark} \label{rem:space-time-white-noise}
We highlight that, the mild solution $\{u(t)\}_{t \in [0, T]}$ of \eqref{eq:intro-SSDVE} driven by space-time
white noise ($Q = I$) can enjoy a positive order of regularity in multiple spatial dimensions ($d > 1$).  To see this,
we first note that Assumption \ref{ass:QWiener} holds in the sense that $\big\|A^{\frac{\gamma-1}{2}} \big\|_{\mathrm{HS}}<\infty $ for $-1\leq \gamma < \frac{2-d}{2}$, $d = 1, 2,3$, by taking the
asymptotics of the eigenvalues of $A$ into account. If $d \in \{ 2, 3\}$, then the estimate
\eqref{eq:thm1-u-spatial2} ensures that $\{u(t)\}_{t \in [0, T]}$ can have a positive order of
spatial regularity since $\gamma + 1 > 0$.  As a comparison, we recall that the mild solutions of
the stochastic heat equation \cite[Corollary 2.5]{yan2005galerkin} and the stochastic wave equation
\cite[Remark 3.2]{kovacs2010finite} subject to the space-time white noise only survive in one spatial dimension.
\end{remark}
\begin{remark}
Here and below,  $C$ denotes a generic positive constant that may vary from line to line,
depending on $T, K$ and $\|A^{\frac{\gamma-1}{2}}Q^{\frac{1}{2}}\|_{\mathrm{HS}}$,
but independent of step-sizes $h, k$. In addition, we make further comments on the initial data.
Since our main interest lies in the influence
due to the presence of the noise,  we work with smooth initial data here and below (e.g., $ \varphi \in L^2 (\Omega;  \dot{H}^{\gamma+1} ) $, $\psi  \in L^2( \Omega; \dot{H}^{\gamma} )$).
However, as indicated in \cite{larsson1991finite},  such conditions can be relaxed with nonsmooth initial data but at the cost of nonuniform error constants $C$ blowing up as $T \rightarrow 0$.
\end{remark}
%

In order to prove Theorem \ref{thm:main-regularity}, we need some properties of the semigroup $\mathcal{S}(t)$, which rely on properties of the corresponding linear deterministic strongly damped wave equation
 \begin{eqnarray}\label{eq:linear-determ-eq}
  \left\{\begin{array}{ll}
u_{tt} + \alpha Au_t + Au = 0,& \;  t \in (0,\; T],\\
  u(0) = u_0, \; u_t ( 0 ) = v_0. &
  \end{array}\right.
  \end{eqnarray}

%
As mentioned earlier, the linear problem has been examined in \cite{larsson1991finite} and some spatial regularity results of the solution are already available there. Nevertheless, they are far from enough for our analysis in this work and we have to develop some new further regularity results. To begin with, we recall the following spatial regularity result from \cite[Lemma 2.3]{larsson1991finite}.
\begin{lemma}\label{lem:S-GROUP-spatial-regularity}
Let $u(t)$ be the solution of the strongly damped wave equation (\ref{eq:linear-determ-eq}). For  any integer $j\geq 0$, real numbers $\rho \in \mathbb{R} $ and $\sigma\in [0,2]$,  we have
\begin{align}
|D_t^ju(t)|_\rho+|D_t^{j+1}u(t)|_{\rho-\sigma}
       &\leq c \, t^{-j}(|u_0|_\rho+|v_0|_{\rho-\sigma}),  \quad \text{ for } \: t>0.
\label{eq:S-GROUP-SP-Regu2}
\end{align}
\end{lemma}
%
%
Furthermore, we need the following integral versions of spatial regularity results.
\begin{lemma}\label{lem:linear-determin-spatial}
Let $u(t)$ be the solution of the strongly damped wave equation (\ref{eq:linear-determ-eq}), then it holds that
\begin{align}
\alpha \int_0^t|u_t(s)|_\beta^2 \, \mathrm{d}s
       \leq& \tfrac12 ( |u_0|_\beta^2 + |v_0|_{\beta-1}^2 ), \quad  \forall \, \beta\in \mathbb{R},
       \label{L-Determin-Spa1}\\
\alpha t  |u_t(t)|_{\beta}^2 + \int_0^ts| u_{tt}(s)|_{\beta-1}^2 \, \mathrm{d} s
       \leq & C (|u_0|_{\beta}^2 + |v_0|_{\beta-1}^2), \quad \forall \, \beta\in \mathbb{R}.
       \label{L-Determin-Spa2}
\end{align}
\end{lemma}
{\it Proof of Lemma \ref{lem:linear-determin-spatial}.}
To prove (\ref{L-Determin-Spa1}),
we multiply both sides of (\ref{eq:linear-determ-eq}) by $A^{ \beta-1 } u_t $ to obtain
\begin{align}
\frac12 \frac{ \text{d} }{ \text{dr} }|u_t(r)|_{\beta-1}^2 + \alpha |u_t(r)|_\beta^2 + \frac12 \frac{ \text{d} }{ \text{dr} }|u(r)|_\beta^2
=0,
\end{align}
which, after integration over $[s, t]$, suggests that
\begin{align}\label{eq:integrand-spatial-rugularity-v}
|u(t)|_\beta^2+|u_t(t)|_{\beta-1}^2+2\alpha\int_s^t|u_t(r)|_\beta^2\,\mathrm{d}r =  |u(s)|_\beta^2+|u_t(s)|_{\beta-1}^2.
\end{align}
Taking $s=0$ implies \eqref{L-Determin-Spa1} straightforwardly.
To validate \eqref{L-Determin-Spa2},  we multiply (\ref{eq:linear-determ-eq}) by $sA^{\beta -1}u_{tt}$  and do some manipulations to arrive at
\begin{equation}
\begin{split}
 s(u_{tt}, A^{\beta-1}u_{tt}) +  \frac{\dd}{\dd s} \big( \tfrac{\alpha}{2}sAu_t + sAu, A^{\beta-1} u_{t} \big)
 =
(s + \tfrac{\alpha}{2})(Au_t, A^{\beta-1} u_{t} )
+
(A u, A^{\beta-1} u_{t}).
\end{split}
\end{equation}
Similarly as before,  by integration over $[0, t]$ and using (\ref{eq:S-GROUP-SP-Regu2})
and (\ref{L-Determin-Spa1}), one can derive that
\begin{align}
\int_0^t s| u_{tt}(s)|_{\beta-1}^2 \, \dd s  + \tfrac{\alpha t }2 |u_t(t)|_{\beta}^2
 &  = - (tu(t),A^{\beta} u_{t}(t))
+\int_0^t \left[ (s+\tfrac{\alpha}{2}) | u_t (s) |^2_{\beta} + ( u, A^{\beta} u_{t}) \right]\, \dd s
\nonumber\\ &  \leq
|u(t)|_{\beta} \, t |u_t(t)|_{\beta} + (T + \tfrac{\alpha+1}{2} ) \int_0^t |u_t(s)|_{\beta}^2 \, \dd s
+
\tfrac12 \int_0^{t}\ |u(s)|^2_{\beta} \, \dd s
\nonumber\\ &  \leq
C ( |u_0|^2_{\beta} + |v_0|^2_{\beta -1} ) + C ( |u_0|_{\beta}^2 + |v_0|^2_{\beta - 1} )
+ C T ( |u_0|_{\beta}^2 + |v_0|^2_{\beta -1} )
\nonumber\\ & \leq
C ( |u_0|_{ \beta }^2 + |v_0|_{ \beta -1 }^2 ).
\end{align}
This thus concludes the proof of this lemma.
$\square$

Based on the above spatial regularity results, one can tackle the temporal regularity properties.
%
\begin{lemma}\label{Lem:timez}
Let $u(t)$ be the solution of the equation (\ref{eq:linear-determ-eq}). For  $0\leq s<t\leq T$,  we have
\begin{align}
\|u(t)-u(s)\|\leq C(t-s)^{\frac \mu 2}(|u_0|_\mu+|v_0|_{\mu-2}),\quad \forall \mu \in [0,2],
\label{eq:linear-determ-u-time-regularity}\\
\|u_t(t)-u_t(s)\|\leq Cs^{-\nu}(t-s)^{\nu}(|u_0|_2+\|v_0\|), \quad \forall \nu\in[0,1],
\label{eq:linear-determ-v-time-regularity}\\
\|u_t(t)-u_t(s)\|\leq C(t-s)^{\frac\mu2}(|u_0|_\mu+|v_0|_\mu), \quad \forall \mu\in[0,2],
\label{eq:linear-determ-v-time-regularity1}\\
\int_s^t\|u_t(r)\|^2 \,\mathrm{d} r \leq C(t-s)^{\nu}(|u_0|^2_\nu+|v_0|^2_{\nu-1}),
\quad \forall \nu\in[0,1].
\label{eq:linear-determ-v-time-regularity-integrad}
\end{align}
\end{lemma}
{ \it  Proof of Lemma \ref{Lem:timez}.}
Thanks to interpolation theory, we only need to verify (\ref{eq:linear-determ-u-time-regularity}) for the two cases $\mu=0$ and $\mu=2$. With the aid of \eqref{eq:S-GROUP-SP-Regu2} with $\rho=j=0, \sigma = 2$, one can see that
\begin{align}
\label{eq:u-diff1}
\|u(t)-u(s)\| \leq \|u(t)\| + \|u(s)\| \leq 2c ( \|u_0 \| + |v_0|_{-2} ).
\end{align}
Likewise,  using (\ref{eq:S-GROUP-SP-Regu2}) with $\rho=\sigma =2, j=0$ leads us to
\begin{align}\label{jiezhe2}
\|u(t)-u(s)\|\leq\int_s^t \|u_t(r)\|\,\dd r  \leq  c |t-s| ( |u_0|_2 + \|v_0 \| ).
\end{align}
%
With regard to (\ref{eq:linear-determ-v-time-regularity}), in the same manner we use \eqref{eq:S-GROUP-SP-Regu2} with $\rho=j=\sigma=0$ to infer that
\begin{equation}\label{eq:v-holder-regularity-2}
\| u_t (t) - u_t (s) \| \leq  \| u_t (t) \| + \| u_t (s) \|
   \leq
       2 c  ( \| u_0\| + \| v_0 \| )
   .
\end{equation}
At the same time, due to  (\ref{eq:S-GROUP-SP-Regu2}) with $\rho=\sigma =2, j=1$,  we get
\begin{align}\label{eq:v-holder-regularity}
\| u_t(t)-u_t(s) \|
    \leq \int_s^t  \| u_{tt} (r) \| \, \dd r
     \leq c \int_s^t  r^{-1}  ( |u_0|_2 + |v_0 |_0 ) \, \dd r
      \leq  c s^{-1} |t-s| ( |u_0|_2 + \|v_0 \| ).
\end{align}
To show (\ref{eq:linear-determ-v-time-regularity1}), we recall $u_{tt} + \alpha Au_t + Au = 0$ and apply (\ref{eq:S-GROUP-SP-Regu2}) with $\rho=2$, $j=\sigma=0$ to obtain
\begin{align}
\| u_t(t)-u_t(s) \|
    \leq \int_s^t  \| u_{tt} (r) \| \, \dd r
    \leq \int_s^t(\alpha | u_t (r) |_2 +|u(r)|_2)\, \dd r\leq C(t-s)(|u_0|_2+|v_0|_2),
\end{align}
which combined with  (\ref{eq:v-holder-regularity-2}) implies (\ref{eq:linear-determ-v-time-regularity1})
by interpolation. Finally, the proof of (\ref{eq:linear-determ-v-time-regularity-integrad}) for  the cases
$\nu=0$ and $\nu=1$ are, respectively, direct consequences of (\ref{eq:integrand-spatial-rugularity-v})
with $\beta = 0$ and  (\ref{eq:S-GROUP-SP-Regu2}) with $j= 0, \rho=\sigma=1$.
$\square$
%

At this stage we are ready to associate the above regularity results with properties of the semigroup. To this end, we come back to the linear problem \eqref{eq:linear-determ-eq} and reformulate it as a system of first order
\begin{eqnarray}\label{eq:linear-system-1order}
\left\{\begin{array}{ll} w_t + \mathcal{A} w=0,\qquad & t \in  (0, T],\\
w(0)=w_0,\;&
\end{array}\right.
\end{eqnarray}
where we denote $v := u_t, w=(u,v)'$ and $w_0=(u_0, v_0)'$.
In terms of  the semigroup $\mathcal{S}(t)$, the solution of  \eqref{eq:linear-system-1order} is given by
$w(t)= (u(t),v(t))' = \mathcal{S}(t)w_0$ for $w_0\in  \dot{H}^{s} \times  \dot{H}^{s-\sigma}, s \in \mathbb{R},
\sigma \in [0, 2]$. For two Hilbert spaces $H_i$, $i=1,2$, we additionally introduce two operators
$P_1$ and $P_2$ defined by
\[
P_ix=x_i,\;\forall x=\left(x_1,x_2\right)'\in H_1\times H_2.
\]
Noting that $u(t) = P_1 \mathcal{S}(t)w_0$ and $v(t) = P_2 \mathcal{S}(t)w_0$, one can reformulate the above regularity results in a semigroup way.
For example, \eqref{eq:S-GROUP-SP-Regu2} with $j = 0$ can be rewritten as
\begin{equation} \label{eq:spatial-regularity-semigroup}
|P_1 \mathcal{S}(t) w_0 |_\rho + | P_2 \mathcal{S}(t) w_0 |_{\rho-\sigma}
       \leq C ( | P_1 w_0 |_\rho + | P_2 w_0 |_{\rho-\sigma}),  \quad  \rho \in \mathbb{R},\, \sigma \in [0, 2],\quad t>0.
\end{equation}
Moreover, Lemmas \ref{lem:linear-determin-spatial}-\ref{Lem:timez} suggest that
\begin{align}
\alpha t  | P_2 \mathcal{S}(t) w_0 |_{\beta}^2
  + \alpha \smallint_0^t | P_2 \mathcal{S}( s ) w_0 |_\beta^2 \, \mathrm{d}s
       \leq C( | P_1 w_0 |_\beta^2 + | P_2 w_0 |_{\beta-1}^2), & \quad  \forall \beta\in \mathbb{R},
       \label{eq:integral-spatial} \\
\| P_1 \mathcal{S}(t) w_0 - P_1 \mathcal{S}(s) w_0 \|\leq C(t-s)^{\frac \mu 2}
(|P_1 w_0|_\mu + | P_2 w_0 |_{\mu-2}), & \quad \forall \mu \in [0,2],
\label{eq:linear-deter-u-time-regularity-S}\\
\| P_2 \mathcal{S}(t) w_0 - P_2 \mathcal{S}(s) w_0 \|\leq  Cs^{-\nu}(t-s)^{\nu}
(| P_1 w_0 |_2 + \| P_2 w_0 \| ), & \quad \forall \nu\in[0,1],
\label{eq:linear-deter-v-time-regularity-S}\\
\| P_2 \mathcal{S}(t) w_0 - P_2 \mathcal{S}(s) w_0 \|\leq  C(t-s)^{\frac\mu2}
(| P_1 w_0 |_\mu + | P_2 w_0 |_{\mu}), & \quad \forall \mu\in[0,2],
\label{eq:linear-deter-v-time-regularity-Ss}\\
\smallint_s^t\|P_2\mathcal{S}(r)w_0\|^2\,\dd r
\leq C(t-s)^\nu (|P_1w_0|^2_\nu+|P_2w_0|^2_{\nu-1}),& \quad \forall \nu\in[0,1].
\label{eq:linear-deter-v-time-regularity-Sss}
\end{align}

{\it Proof of Theorem \ref{thm:main-regularity}.}
Let $\varrho \in [-1, \gamma]$ and let $H=\dot{H}^{\varrho + 1} \times \dot{H}^{\varrho -1}$ be equipped with the
norm $\|X\|_H^2 =  | u |_{\varrho + 1}^2 + | v |_{\varrho -1}^2,\forall X=(u,v)'\in H$. It is easy to check that $H$ is a
separable Hilbert space. Then we show existence of a unique mild solution in $H$.
Since Assumption \ref{ass:F} holds with $\gamma \in [-1, 1]$ and due to the definition of $\|\cdot\|_H$, we realize that, for any $X=(u,v)'$ and $X_i=(u_i,v_i)', i= 1, 2$,
\begin{align}
\| \mathbf{F} ( X_1 ) - \mathbf{F} ( X_2 ) \|_H & = | F(u_1)-F(u_2) |_{\varrho-1}  \leq C\|u_1-u_2\| \leq C\|X_1-X_2\|_H,
\label{eq:F-condition1}\\
\label{eq:F-condition2}
\| \mathbf{F} ( X ) \|_{H} &= | F(u) |_{\varrho - 1}
\leq C ( 1 + \|u\| )
\leq C (1 + \|X\|_H).
\end{align}
Additionally, the definition of the Hilbert-Schmidt norm and Assumption \ref{ass:QWiener}
enable us to deduce that
\begin{equation}
\begin{split}
\| \mathbf{B}\|^2_{\mathcal{L}_0^2}
& =
\sum_{i \in \mathbb{N} }
\| A^{\frac{\varrho - 1}{2} } Q^{\frac12} e_i \|^2
=
  \| A^{\frac{\varrho - 1}{2}}Q^{\frac{1}{2}}\|_{\mathrm{HS}}^2
\leq
 \| A^{\frac{\gamma - 1}{2}}Q^{\frac{1}{2}}\|_{\mathrm{HS}}^2
 < \infty.
\label{tiaojian3}
\end{split}
\end{equation}
%
%
In view of Theorem 7.4 in \cite{da2014stochastic}, (\ref{eq:F-condition1})-(\ref{tiaojian3})
together with the fact that $\mathcal{A}$ generates an analytic semigroup $\mathcal{S}(t)$ in $H$ guarantee  a unique mild solution given by (\ref{eq:mild-SPDE}), which satisfies
\begin{align}\label{eq:mild-spatial-regularity-first}
\sup_{t\in[0,T]}
\Big(
 \|u(t)\|_{L^2(\Omega;\dot{H}^{\varrho + 1} )} ^ 2 + \| u_t(t) \|_{L^2(\Omega;\dot{H}^{ \varrho-1})}^2
 \Big)
\leq
C \Big( 1 + \|\varphi\|^2_{ L^2(\Omega; \dot{H}^{\varrho + 1} ) }
+
\| \psi \|^2_{ L^2(\Omega; \dot{H}^{\varrho-1} ) } \Big),
\end{align}
for all $\varrho \in [-1, \gamma]$. Taking $\varrho = \gamma$  thus confirms \eqref{eq:thm1-u-spatial2}.
As another consequence of \eqref{eq:mild-spatial-regularity-first}, we have
\begin{equation} \label{eq:F-moment-bound}
\sup_{s \in [0, T] }
\| F(u(s)) \|_{L^2(\Omega;\dot{H}^0)}
\leq
C \big( 1 + \|\varphi\|_{L^2(\Omega; \dot{H}^{0})}
 + \| \psi \|_{L^2(\Omega; \dot{H}^{-2})} \big).
\end{equation}
Concerning the temporal regularity, we apply the It\^{o} isometry to obtain for all $\varrho\in[-1,\gamma]$
\begin{align}\label{eq:u(t)-u(s)}
\| u(t)-u(s) \|_{L^2(\Omega;\dot{H}^0)}
    \leq &
         \|P_1( \mathcal{S}(t-s) - \mathcal{S}(0) ) X(s) \|_{L^2( \Omega;\dot{H}^0 )}
         +
         \int_s^t \|P_1\mathcal{S}(t-r)\mathbf{F}(X(r))\|_{L^2(\Omega;\dot{H}^0)}\,\dd r
     \nonumber \\ &
           +
            \left(\int_s^t \| P_1 \mathcal{S}(t-r) \mathbf{B} Q^{\frac{1}{2}}\|_{\mathrm{HS}}^2 \,\dd r\right)^{\frac12}
     := \mathbb{I}_1 + \mathbb{I}_2 + \mathbb{I}_3.
\end{align}
Combining \eqref{eq:linear-deter-u-time-regularity-S} and  \eqref{eq:mild-spatial-regularity-first} shows
\begin{align}
\mathbb{I}_1 \leq &
          C (t-s)^{\frac{\varrho+1}{2} } \big( \| P_1 X(s) \|_{L^2(\Omega; \dot{H}^{\varrho+1} )}
             + \| P_2 X(s) \|_{L^2 (\Omega; \dot{H}^{\varrho-1} ) } \big)
       \nonumber \\ \leq &
           C \big( 1 + \|\varphi\|_{L^2(\Omega; \dot{H}^{\varrho + 1})}
             + \|\psi\|_{L^2(\Omega; \dot{H}^{\varrho -1})}  \big)
                (t-s)^{\frac{\varrho+1}{2} }.
\end{align}
For the estimate of $\mathbb{I}_2$, one can recall \eqref{eq:spatial-regularity-semigroup} with $\rho = \sigma = 0$, together with \eqref{eq:F-moment-bound} to derive
\begin{equation}
\mathbb{I}_2 \leq C\int_s^t \| F(u(r)) \|_{L^2(\Omega;\dot{H}^0)}\,\dd r
\leq
C \big( 1 + \|\varphi\|_{L^2(\Omega; \dot{H}^{0})}
             + \|\psi\|_{L^2(\Omega; \dot{H}^{-2})}  \big) (t-s).
\end{equation}
To treat the remaining term $\mathbb{I}_3$, we first use the definition of the Hilbert-Schmidt norm,
\eqref{eq:linear-deter-u-time-regularity-S} to get
\begin{equation}
\| P_1 \mathcal{S}(t-r) \mathbf{B} Q^{\frac{1}{2}}\|_{\mathrm{HS}}
=
\| P_1(\mathcal{S}(0) -\mathcal{S}(t-r)) \mathbf{B} Q^{\frac{1}{2}}\|_{\mathrm{HS}}
\leq
C  ( t - r )^{ \frac{\varrho + 1}{2} } \| A^{\frac {\varrho -1}{2} }
                Q^{\frac{1}{2}} \|_{\mathrm{HS}},
\end{equation}
where we also used the fact $P_1 \mathcal{S}(0) \mathbf{B} Q^{\frac12} \phi  = P_1 \mathbf{B} Q^{\frac12} \phi  = 0$ for $\phi \in \dot{H}^0 $. This yields the estimate of $\mathbb{I}_3$:
\begin{equation}
\mathbb{I}_3 \leq C  ( t - r )^{ \frac{\varrho + 2}{2} } \| A^{\frac {\varrho -1}{2} }
                Q^{\frac{1}{2}} \|_{\mathrm{HS}}.
\end{equation}
Putting the above three estimates together implies
\begin{equation}\label{eq:u(t)-u(s).final}
\| u(t)-u(s) \|_{L^2(\Omega;\dot{H}^0)}
         \leq
             C \big( 1 + \|\varphi\|_{L^2(\Omega; \dot{H}^{\varrho + 1})}
              + \|\psi\|_{L^2(\Omega; \dot{H}^{\varrho -1})}  \big)
                (t-s)^{\frac{\varrho+1}{2} },
\end{equation}
which gives \eqref{eq:thm1-u-temporal2} by taking $\varrho = \gamma$.
%
Next, we shall look at the regularity of $\{u_t (t) \}_{t \in [0, T]}$ when Assumption \ref{ass:QWiener}
holds for  $\gamma \in [0, 1]$. The It\^{o} isometry ensures
\begin{equation}
\begin{split}
\| u_t(t) \|_{L^2(\Omega;\dot{H}^{\gamma})} & \leq
\|P_2\mathcal{S}(t)X_0\|_{L^2(\Omega;\dot{H}^{\gamma})}
+
\int_0^t \big \| P_2\mathcal{S}(t-s)\mathbf{F}(X(s)) \big\|_{L^2(\Omega;\dot{H}^\gamma)}\, \dd s
\\ & \quad +
\left(\int_0^t \big\| A^{\frac{\gamma}{2}}P_2\mathcal{S}(t-s)\mathbf{B} Q^{\frac12}
\big\|^2_{ \text{HS} }\, \dd s \right)^{\frac{1}{2}}
:= I_1+I_2+I_3.
\end{split}
\end{equation}
In what follows, we estimate $I_1, I_2, I_3$ separately. The use of \eqref{eq:spatial-regularity-semigroup}
with $\rho= \gamma, \sigma = 0$ guarantees that
\begin{align} \label{eq:v-spatial-regul-I1}
I_1 \leq C(\|\varphi\|_{L^2(\Omega; \dot{H}^{\gamma})}+\|\psi\|_{L^2(\Omega; \dot{H}^{\gamma})}).
\end{align}
Considering \eqref{eq:integral-spatial} with $\beta = \gamma$ and \eqref{eq:F-moment-bound} shows that
\begin{align}
I_2
\leq C\int_0^t(t-s)^{-\frac 1 2} \|A^\frac{\gamma-1}{2}F(u(s))\|_{L^2(\Omega;\dot{H}^0)}\,\dd s
\leq
C \big( 1 + \|\varphi\|_{L^2(\Omega; \dot{H}^{0})}
 + \| \psi \|_{L^2(\Omega; \dot{H}^{-2})} \big).
 \label{eq:v-spatial-regul-I2}
\end{align}
Finally, using  \eqref{eq:integral-spatial} with $\beta = \gamma$ and Assumption \ref{ass:QWiener} yields
\begin{align}
I_3
  \leq C\|A^{\frac{\gamma-1}{2}}Q^{\frac{1}{2}}\|_{\mathrm{HS}}
< \infty.
\end{align}
This together with \eqref{eq:v-spatial-regul-I1} and \eqref{eq:v-spatial-regul-I2} gives \eqref{eq:thm1-v-spatial2}.
%
To prove (\ref{eq:thm1-v-temporal2}),
 it holds  by (\ref{eq:spatial-regularity-semigroup}),  (\ref{eq:linear-deter-v-time-regularity-Ss}) and (\ref{eq:linear-deter-v-time-regularity-Sss})
\begin{align}
\| u_t(t)-u_t(s) \|_{L^2(\Omega;\dot{H}^0)}
    \leq &
         \|P_2( \mathcal{S}(t-s) - \mathcal{S}(0) ) X(s) \|_{L^2( \Omega;\dot{H}^0 )}
         +
         \int_s^t \|P_2\mathcal{S}(t-r)\mathbf{F}(X(r))\|_{L^2(\Omega;\dot{H}^0)}\,\dd r
     \nonumber \\ &
           +
         \left(\int_s^t \| P_2\mathcal{S}(t-r) \mathbf{B} Q^{\frac{1}{2}}\|_{\mathrm{HS}}^2 \,\dd r\right)^{\frac12}
      \nonumber \\  \leq &
          C (t-s)^{\frac{\gamma}{2} } \big( \| P_1 X(s) \|_{L^2(\Omega; \dot{H}^{\gamma} )}
             + \| P_2 X(s) \|_{L^2 (\Omega; \dot{H}^{ \gamma} ) } \big)
      \nonumber \\ & +
          C \int_s^t \|F(u(r)) \|_{L^2(\Omega;\dot{H}^0)}\,\dd r
            + C ( t - s )^{ \frac\gamma2  } \| A^{\frac {\gamma -1}{2} }
                Q^{\frac{1}{2}} \|_{\mathrm{HS}}
       \nonumber \\ \leq &
           C \big( 1 + \|\varphi\|_{L^2(\Omega; \dot{H}^{\gamma })}
             + \|\psi\|_{L^2(\Omega; \dot{H}^{\gamma })}  \big)
                (t-s)^{\frac{ \gamma}{2} },
\end{align}
where we also used the fact $P_1 \mathbf{B} Q^{\frac12} \phi  = 0$, (\ref{eq:thm1-v-spatial2}),
(\ref{eq:F-moment-bound}) and (\ref{eq:mild-spatial-regularity-first}) with $\varrho=\gamma-1$, $\gamma\in [0,1]$.
$\square$


\section{Error estimates for the finite element semi-discretization and full-discretization of the deterministic linear problem}
\label{sec:Error-estimate-determin}
In this section, we  consider  the semi-discrete and full-discrete finite element approximations of  the deterministic linear strongly damped wave equation (\ref{eq:linear-determ-eq}).  A variety of error estimates will be derived, which play an important role in the mean-square convergence analysis of the finite element method for the stochastic strongly damped wave equation. 

For simplicity of presentation, we assume that $ L = \Delta$ in the following. Let $V_h\subset H_0^1(\mathcal{D}), \, h \in (0, 1]$ be the space of continuous functions that are piecewise linear over the triangulation $\mathcal{T}_h$ of $\mathcal{D}$.  Then 
we define the discrete Laplace operator $A_h  \colon V_h \rightarrow V_h $  by
\begin{equation}\label{eq:defn-Ah}
(A_h v_h,\chi_h)=a(v_h,\chi_h):=(\nabla v_h,\nabla \chi_h),\;\forall \ v_h,\chi_h\in V_h.
\end{equation}
Additionally we introduce a Riesz representation operator $\mathcal{R}_h: H_0^1 (\mathcal{D}) \rightarrow V_h$ defined by
\begin{equation}\label{asr5}
a(\mathcal{R}_h v,\chi_h)=a(v,\chi_h),\;\forall \, v\in H_0^1(\mathcal{D}),\;\forall \,\chi_h\in V_h,
\end{equation}
and a generalized projection operator $\mathcal{P}_h:\dot{H}^{-1}\rightarrow V_h$ given by
\begin{align}
(\mathcal{P}_h v, \chi_h)=(v, \chi_h),\;\forall \, v\in \dot{H}^{-1},\; \forall \, \chi_h\in V_h.
\end{align}
It is well-known that (see e.g., (2.15)-(2.16) in \cite{andersson2012weak}) the operators $\mathcal{P}_h$
and $\mathcal{R}_h$ defined as above satisfy
 \begin{align}
\big\| A^{\frac s 2}(I-\mathcal{R}_h)A^{-\frac r 2} \big\|_{ \mathcal{L}( \dot{H}^0 ) } \leq & C h^{r-s},\quad 0\leq s \leq 1 \leq r\leq 2,\label{IIII1}\\
\big\| A^{\frac s 2}(I-\mathcal{P}_h)A^{-\frac r 2} \big\|_{ \mathcal{L}( \dot{H}^0 ) } \leq & C h^{r-s}, \quad 0\leq s \leq 1,\quad   0\leq s  \leq r \leq 2.
\label{IIII3}
\end{align}
Moreover, the operators $A$ and $A_h$ obey
  \begin{equation}\label{asr3}
C_1 \| A_h^{\frac \gamma 2}\mathcal{P}_hv \| \leq \| A^{\frac {\gamma} {2}}v \|
\leq C_2  \|A_h^{\frac \gamma 2}\mathcal{P}_hv \| ,\quad  v\in \dot{H}^{-\gamma},\; \gamma\in[-1,1].
 \end{equation}
Furthermore, we denote by $T \colon L_{2}( \mathcal{D} ) \rightarrow H^{2}(\mathcal{D}) \cap H^{1}_{0}(\mathcal{D})$ the solution operator of the equation $Au=f$ and $T_{h}:L_{2}( \mathcal{D} ) \rightarrow V_h$ approximation of $T$, so that
\begin{equation}
a(T_hf,\chi_h)=(f,\chi_h),\quad \forall \, f \in L_{2}( \mathcal{D} ),  \chi_h \in V_{h}.
\end{equation}
By the definition of the operator  $\mathcal{R}_h$, we observe $T_h= A_h^{-1} \mathcal{P}_h =\mathcal{R}_hT$,
and $T_h$ is self-adjoint, positive semi-definite on $L_2( \mathcal{D} )$, and positive definite on $V_h$.
Furthermore, as a consequence of \eqref{IIII1} we have
\begin{align}\label{dadlal}
\|(T_h-T)f\|\leq Ch^{s}\|f\|_{s-2}, \quad f\in \dot{H}^{s-2},\;  s \in [1, 2].
\end{align}

\subsection{Error estimates of  semidiscrete scheme}
In this subsection, we focus on the semi-discrete finite element approximation of the deterministic linear problem \eqref{eq:linear-determ-eq} and prove some useful estimates. We mention that such error estimates for the semi-discrete scheme will be derived based on energy arguments and some known results in \cite{larsson1991finite}.

Note first that the weak variational form of \eqref{eq:linear-determ-eq} is to find $(u,v)'\in H_0^1( \mathcal{D} ) \times H_0^1( \mathcal{D} ) $ such that
\begin{eqnarray}\label{EA222}\label{sspdr4}
\left\{\begin{array}{ll}(v_t, \vartheta_1)+\alpha(\nabla v,\nabla \vartheta_1)+(\nabla u,\nabla \vartheta_1)=0,&\;\vartheta_1\in H_0^1( \mathcal{D} ) ,\\
(\nabla u_t,\nabla\vartheta_2)-(\nabla v,\nabla\vartheta_2)=0,&\;\vartheta_2\in H_0^1 ( \mathcal{D} ).
\end{array}\right.
\end{eqnarray}
The corresponding semidiscrete finite element method is thus to find $(u_h,v_h)'\in V_h\times V_h$
such that
\begin{eqnarray}\label{EA222}\label{eq:FEM-weak}
\left\{\begin{array}{ll}(v_{h,t}, \vartheta_1)+\alpha(\nabla v_h,\nabla \vartheta_1)+(\nabla u_h,\nabla \vartheta_1)=0,&\;\vartheta_1\in V_h,\\
(\nabla u_{h,t},\nabla\vartheta_2)-(\nabla v_h,\nabla\vartheta_2)=0,&\;\vartheta_2\in V_h.
\end{array}\right.
\end{eqnarray}
In terms of the discrete Laplace operator defined by \eqref{eq:defn-Ah}, we can equivalently write \eqref{eq:FEM-weak} as
\begin{eqnarray}\label{eq:FEM-system}
\left[\begin{array}{c}u_{h,t}(t)\\v_{h,t}(t)\end{array}\right]+
\left[\begin{array}{c c}0 &-I \\A_h
&\alpha A_h\end{array}\right]\biggl[\begin{array}{c}u_h(t)\\v_h(t)\end{array}\biggr]
=0\quad \text{ with }\; u_h(0)=\mathcal{P}_hu_0, \;\; v_h(0) = \mathcal{P}_hv_0.
\end{eqnarray}
Similarly as \eqref{eq:linear-system-1order}, we can also reformulate it as
\begin{align}\label{eq:FEM-system1}
\left\{\begin{array}{ll}
w_{h,t}+\mathcal{A}_hw_h=0,&\;t \in (0, T], \\
w_h(0)=w_{0h},&
\end{array}\right.
\end{align}
where we denote $w_h(t)=(u_h(t),v_{h}(t))'\in V_h\times V_h,
w_{0h}=(\mathcal{P}_hu_0,\mathcal{P}_hv_0)'$ and
\[\mathcal{A}_h:=\biggl[\begin{array}{c
c}0&-I\\A_h&\alpha A_h\end{array}\biggr].\]
Here $-\mathcal{A}_h$ generates an analytic semigroups $\mathcal{S}_h(t)$  in $V_h\times V_h$ supplied with the norm of $L_2(\mathcal{D}) \times L_2(\mathcal{D})$ \cite{larsson1991finite}.
Let $\mathbf{P}_h$ denote a projection operator from $\dot{H}^{-1}\times \dot{H}^{-1}$ to $V_h\times V_h$ defined by $\mathbf{P}_hx=(\mathcal{P}_hx_1,\mathcal{P}_hx_2)'$, $\forall x=(x_1,x_2)\in \dot{H}^{-1}\times \dot{H}^{-1}$. Then the solution of (\ref{eq:FEM-system1}) can be written as $w_h(t)=\mathcal{S}_h(t)\mathbf{P}_hw_0$.
The following results can be regarded as an extension of error estimates
of integer order (i.e., $\beta = 2, q =2$) in \cite{larsson1991finite} to cover intermediate cases.
%
%
%
\begin{theorem}\label{thm:error-FEM1}
Let $w_h(t)=(u_h(t),v_h(t))'$ and $w=(u(t),v(t))'$ be the solutions of  (\ref{eq:FEM-system}) and (\ref{eq:linear-system-1order}), respectively. Let the setting in the beginning of section \ref{sec:Error-estimate-determin} be fulfilled and define an error operator as
\begin{equation}
F_h(t)w_0:=w_h(t)-w(t)=(\mathcal{S}_h(t)\mathbf{P}_h - \mathcal{S}(t))w_0, \quad w_0=(u_0,v_0)'.
\end{equation}
Then it holds that
\begin{align}
\|P_1F_h(t)w_0\|\leq Ch^\beta(|u_0|_\beta + |v_0|_{\beta-2}), & \quad 1\leq \beta\leq 2,
\label{eq:error-FEM-u}\\
\|P_2F_h(t)w_0\|\leq C h^q(|u_0|_q+|v_0|_q), & \quad  0\leq q\leq 2
\label{eq:error-FEM-v}.
\end{align}
\end{theorem}
{\it Proof of  Theorem \ref{thm:error-FEM1}.}
%
The estimate (\ref{eq:error-FEM-u}) for the special case $\beta = 2$ can be found in \cite[Theorem 3.1]{larsson1991finite}.  Following the basic lines there and taking \eqref{dadlal} into account, one can readily
justify (\ref{eq:error-FEM-u}) for the intermediate cases $\beta \in [1,2]$.
For (\ref{eq:error-FEM-v}), the case $q=0$ and $q=2$ can be
immediately achieved by applying the stability of $\mathcal{S}(t)$ and $\mathcal{S}_h(t) \mathbf{P}_h$ in $L_2(\mathcal{D})
\times L_2(\mathcal{D}) $, and  Theorem 3.4 in \cite{larsson1991finite}, respectively.
The interpolation theory thus results in \eqref{eq:error-FEM-v} for the general case.
$\square$
\begin{remark}\label{remark:banqdr2}
It is worthwhile to point out that (\ref{eq:error-FEM-u}) can not hold for the range $\beta \in [0, 1)$
because \eqref{dadlal} is valid only for $s \in [1, 2]$.
Indeed, the projection $\mathbf{P}_h w_0$ is not well-defined if $v_0 \in \dot{H}^{\beta -2}$, $\beta <1$.
Also, this comment applies to the error estimate \eqref{eq:error-full-u} below for the full discretization.
\end{remark}

Subsequently, we will present a completely new error estimate of integral form, which requires weaker regularity assumption
on $v_0$ and can not be derived directly from existing results in \cite{larsson1991finite}.
Moreover, a nonsmooth data error estimate like \eqref{eq:error-FEM-v2} is obtained.
\begin{theorem}\label{thm:error-FEM2}
 Under the assumptions stated in Theorem \ref{thm:error-FEM1},  it holds that
\begin{align}
\|P_1F_h(t)w_0\| +
\left( \int_0^t \|P_2F_h(s)w_0\|^2 \,\mathrm{d}s \right)^{\frac12}
\leq
Ch^\beta(|u_0|_\beta+|v_0|_{\beta-1}),& \quad  0\leq \beta\leq 2,
\label{eq:error-FEM-u-v} \\
\|P_2F_h(t)w_0\| \leq Ch^qt^{-\frac{q-s}{2}}(|u_0|_{s+2}+|v_0|_s), &\quad  0\leq s\leq q\leq 2.
\label{eq:error-FEM-v2}
\end{align}
\end{theorem}
{\it Proof of Theorem \ref{thm:error-FEM2}.}
By interpolation, we only need to show \eqref{eq:error-FEM-u-v} for $\beta=0$ and $\beta=2$.
For the case $\beta=0$, we set $\vartheta_1=T_h v_h$ and $\vartheta_2=T_hu_h$ in  (\ref{EA222})
and add the resulting two equations to get
\begin{equation}
(v_{h, t}, T_h v_h ) + \alpha ( v_h, v_h) + ( u_{h, t}, u_h ) = 0,
\end{equation}
where the definitions of $T_h$ and $A_h$ were also used. Equivalently, we can recast it as
\begin{equation}\label{eq:spr5}
\frac12 \frac{ \text{d} }{ \text{ds} }\|u_h(s)\|^2 + \frac12 \frac{ \text{d} } { \text{ds} }\|T^{\frac{1}{2}}_hv_h(s)\|^2
+\alpha\|v_h(s)\|^2 = 0,
\end{equation}
which, after integration over $[0, t]$ and employing (\ref{asr3}), leads to
\begin{align}\label{sht}
\|u_h(t)\|^2+\|T_h^{\frac{1}{2}}v_h(t)\|^2 + 2\alpha\int_0^t\|v_h(s)\|^2\,\mathrm{d}s
& =
\|u_h(0)\|^2+\|T_h^{\frac{1}{2}}v_h(0)\|^2
\nonumber\\
& =\|\mathcal{P}_hu_0\|^2+\|T_h^{\frac{1}{2}}\mathcal{P}_hv_0\|^2\leq C(\|u_0\|^2+|v_0|_{-1}^2).
\end{align}
In the same spirit as \eqref{sht}, one can derive that
\begin{align}\label{sht1}
\|u(t)\|^2+\|T^{\frac{1}{2}}v(t)\|^2+2\alpha\int_0^t\|v(s)\|^2\,\mathrm{d}s =  \|u_0\|^2+|v_0|_{-1}^2.
\end{align}
Therefore, combining  \eqref{asr3}, \eqref{sht} and  \eqref{sht1} enables us to get
\begin{align}
\|u(t)-u_h(t)\|^2 + \|v(t)-v_h(t)\|_{-1}^2 + \int_0^t\|v-v_h\|^2\,\mathrm{d}s
\leq C(\|u_0\|^2+|v_0|_{-1}^2).
\end{align}
This verifies \eqref{eq:error-FEM-u-v} in the case $\beta=0$.
To prove (\ref{eq:error-FEM-u-v}) for $ \beta=2 $,
we introduce notations as follows,
\begin{equation} \label{eq:error-decomposition}
\begin{split}
\theta_1:=  u_h - \mathcal{R}_h u,\quad & \rho_1 = (\mathcal{R}_h-I) u,\\
\theta_2:=  v_h - \mathcal{P}_hv,\quad  & \rho_2 = (\mathcal{P}_h-I)v,\\
e_1:= u_h - u =  \theta_1 + \rho_1, \quad  & e_2:=  v_h - v = \theta_2 + \rho_2.
\end{split}
\end{equation}
Subtracting  (\ref{sspdr4}) from $(\ref{EA222})$ shows
\begin{align}
(e_{2,t},\vartheta_1)+\alpha(\nabla e_2, \nabla
\vartheta_1)+(\nabla e_1, \nabla
\vartheta_1)=&0,\quad\forall\vartheta_1\in V_h,\\
(\nabla e_{1,t},\nabla\vartheta_2)-(\nabla e_2,\nabla \vartheta_2)=&0,\quad \forall\vartheta_2\in
V_h.
\end{align}
Inserting $e_i=\theta_i+\rho_i,i=1,2$ thus gives
\begin{align}
(\theta_{2,t}, \vartheta_1)+\alpha(\nabla\theta_2,\nabla
\vartheta_1)+(\nabla\theta_1,\nabla
\vartheta_1)=-( \rho_{2,t},  \vartheta_1)&-\alpha(\nabla\rho_2,\nabla
\vartheta_1)-(\nabla\rho_1,\nabla
\vartheta_1),\;\forall \vartheta_1\in V_h,\\
(\nabla\theta_{1,t},\nabla\vartheta_2)-(\nabla \theta_2,\nabla \vartheta_2)=&-(\nabla\rho_{1,t},\nabla\vartheta_2)+(\nabla \rho_2,\nabla \vartheta_2),\;\; \forall
\vartheta_2\in V_h.
\end{align}
Further,  the orthonormal properties of the operators $\mathcal{R}_h$ and $\mathcal{P}_h$
help us to arrive at
\begin{align}
(\theta_{2,t}, \vartheta_1)+\alpha(\nabla\theta_2,\nabla
\vartheta_1) + (\nabla\theta_1,\nabla
\vartheta_1) = & -\alpha(\nabla \mathcal{R}_h \rho_2,\nabla
\vartheta_1),\quad \forall \vartheta_1\in V_h,\\
(\nabla\theta_{1,t},\nabla\vartheta_2) - (\nabla \theta_2,\nabla \vartheta_2)
= & (\nabla \mathcal{R}_h \rho_2, \nabla \vartheta_2),
\quad\quad\quad \forall \vartheta_2\in V_h.
\end{align}
Setting  $\vartheta_1=T_h \theta_2$,  $\vartheta_2=T_h\theta_1$, adding the resulting two equations,
and taking the definitions of $T_h, A_h$ into account yield that
\begin{equation} \label{eq:semi-discrete-error1}
(\theta_{1,t}, \theta_1)+(\theta_{2,t}, T_h \theta_2)+\alpha(\theta_2,
\theta_2)=-\alpha(\mathcal{R}_h\rho_2,  \theta_2)+(\mathcal{R}_h\rho_2, \theta_1).
\end{equation}
Exploiting similar arguments as before, in conjunction with Cauchy-Schwarz inequality shows
\begin{equation} \label{spr7}
\begin{split}
& \|\theta_1(t)\|^2 +\|T_h^{\frac{1}{2}}\theta_2(t)\|^2 + \alpha\int_0^t\|\theta_2(s)\|^2\, \mathrm{d}s
\\
& \quad \leq (\alpha+1)\int_0^t\|\mathcal{R}_h\rho_2(s)\|^2\,\textmd{d}s + \int_0^t\|\theta_1(s)\|^2\,\textmd{d}s
+  \|\theta_1(0)\|^2,
\end{split}
\end{equation}
where we also used the fact that $\theta_2(0)=v_h(0)-\mathcal{P}_hv(0)=0$.
Applying Gronwall's inequality yields
\begin{equation}
\begin{split}
\|\theta_1(t)\|^2 +\|T_h^{\frac{1}{2}}\theta_2(t)\|^2 &+ \alpha \int_0^t\|\theta_2(s)\|^2\,\mathrm{d}s
\leq
C \left( \int_0^t\|\mathcal{R}_h \rho_2 (s) \|^2\mathrm{d}s + \|\theta_1(0)\|^2 \right)\\
&  \leq Ch^4 \left(\int_0^t|u_t(s)|_2^2\,\mathrm{d}s+|u_0|_2^2\right)
 \leq Ch^4(|u_0|_2^2 + |v_0|_1^2 ),
\end{split}
\end{equation}
where we also used \eqref{L-Determin-Spa1} and the facts that $\|\theta_1(0)\| =
\|\mathcal{P}_h ( \mathcal{R}_h-I)u_0\| \leq Ch^2|u_0|_2$ and that
\begin{equation}
\|\mathcal{R}_h\rho_2\| = \|\mathcal{R}_h (\mathcal{P}_h-I)v \| = \| \mathcal{P}_h (\mathcal{R}_h - I ) v \|
\leq  Ch^2|v|_2.
\end{equation}
Furthermore, \eqref{eq:S-GROUP-SP-Regu2}, \eqref{L-Determin-Spa1}, \eqref{IIII1} and \eqref{IIII3} promise that
\begin{equation}
\| \rho_1 (t) \|^2 + \int_0^t \| \rho_2(s) \|^2 \, \mathrm{d}s  \leq
Ch^4 | u(t) |_2^2 + C h^4 \int_0^t | v(s) |_2^2 \, \mathrm{d}s
\leq  Ch^4 ( |u_0|_2^2 + |v_0|_1^2 ).
\end{equation}
The triangle inequality shows \eqref{eq:error-FEM-u-v} for $\beta=2$ and the interpolation argument finally concludes the proof of \eqref{eq:error-FEM-u-v}.
The assertion  \eqref{eq:error-FEM-v2} can be also deduced by interpolation between $s = 0$ and $s = q$
(see, e.g., \cite[Theorem 3.5]{thomee2006galerkin}). The case $s= q$  is an immediate consequence of
\eqref{eq:error-FEM-v}. To prove \eqref{eq:error-FEM-v2} for $s = 0$, we again use interpolation arguments.
The case $s =0,  q = 0$ is a special case of \eqref{eq:error-FEM-v} and the case $s = 0, q = 2$ is covered by \cite[Theorem 3.1]{larsson1991finite}.
$\square$
%
%

\subsection{Error estimates of the full-discrete scheme}
In this part, we turn our attention to the full-discrete finite element approximation of (\ref{eq:linear-determ-eq}) with
the backward Euler scheme for the time discretization.
With the help of previous findings for the spatially discrete scheme,
we will derive some error estimates for the full-discrete scheme.
We mention that the error analysis in the full-discrete setting is more involved and one
needs to explore further regularity results of the linear deterministic problem
and the analysis relies heavily on energy arguments.

Let $k$ be a time step-size such that $k = \tfrac{T}{N}$, $N \in \mathbb{N}$ and denote
$t_n=nk$, for $0 \leq n\leq N$.  Applying the backward Euler scheme to the
semi-discretization problem \eqref{eq:FEM-system}  gives the full-discrete finite element approximation. More accurately, we are to find $W^n=(U^n, V^n)' \in V_h\times V_h$ such that
\begin{eqnarray}\label{dspr}
\left[\begin{array}{c}U^n\\V^n\end{array}\right]
-
\left[\begin{array}{c}U^{n-1}\\V^{n-1}\end{array}\right]+k\biggl[\begin{array}{c
c}0
&
-
I\\A_h&\alpha A_h\end{array}\biggr]\biggl[\begin{array}{c}{U^n}\\{V^n}\end{array}\biggr]
=0 \quad \text{ with }  \:
\left[\begin{array}{c}U^0 \\ V^0 \end{array}\right] = \left[\begin{array}{c} \mathcal{P}_hu_0
\\ \mathcal{P}_hv_0 \end{array}\right],
\end{eqnarray}
or in a compact way,
\begin{eqnarray}\label{DPf2}
W^n-W^{n-1}+k\mathcal{A}_hW^n=0,\quad W^0 = \mathbf{P}_hw_0.
\end{eqnarray}
Denoting $r(\lambda)=(1+\lambda)^{-1}$, we can rewrite $(\ref{DPf2})$ in
the form
\begin{equation}\label{DPf3}
W^n=r(k\mathcal{A}_h)^n \mathbf{P}_hw_0,
\end{equation}
where the operator $r(k\mathcal{A}_h)^n$ is stable in the following sense  \cite[Lemma 5.2]{larsson1991finite}:
\begin{equation}
\label{DDPf5}
\|r(k\mathcal{A}_h)^n\|_{\mathcal{L}( \dot{H}^0 \times \dot{H}^0 ) } \leq C,\; \text{ for } \; n \geq 0.
\end{equation}
Our aim is thus to analyze various error estimates of $W^n-w(t_n)$.
As the first part, we can derive the time discrete analogue of Theorem \ref{thm:error-FEM1}.
\begin{theorem}\label{thm:error-FEM3}
Let all the conditions in Theorem \ref{thm:error-FEM1} be fulfilled and define
\begin{equation}
\label{eq:error-fem-representation}
F_{kh}^nw_0:= W^n - w(t_n) = (r(k\mathcal{A}_h)^n\mathbf{P}_h-\mathcal{S}(t_n))w_0, \quad w_0=(u_0,v_0)',
\end{equation}
where $W^n$ and $w(t_n)$ obey (\ref{DPf3}) and (\ref{eq:linear-system-1order}), respectively.
Then it holds that
\begin{align}
\|P_1F_{kh}^nw_0\|&\leq C(h^\beta+k^{\frac{\beta}{2}}) (|u_0|_\beta+|v_0|_{\beta-2}), \quad 1
\leq \beta\leq 2,
\label{eq:error-full-u}\\
\|P_2F_{kh}^nw_0\|&\leq C(h^q+k^{\frac{q}{2}})(|u_0|_{q}+|v_0|_q), \quad 0\leq q\leq 2.
\label{eq:error-full-v}
\end{align}
\end{theorem}
%
{\it Proof of Theorem \ref{thm:error-FEM3}.}
Owing to interpolation arguments, we only need to show \eqref{eq:error-full-u} for $\beta=1$ and $\beta=2$.
The latter case $\beta=2$ can be found in \cite[Theorem 5.3]{larsson1991finite}.
For the case $\beta=1$, with the spatial approximation error \eqref{eq:error-FEM-u} at disposal,  it remains to prove
\begin{equation} \label{eq:thm-U-full-estimate1}
\| P_1 F_n \mathbf{P}_h w_0 \|
=
\| P_1 (r(k\mathcal{A}_h)^n - \mathcal{S}_h(t_n)) \mathbf{P}_h w_0 \|
\leq
C k^{\frac12} ( |u_0|_1 + |v_0|_{-1} ),
\end{equation}
where we introduced an error operator $F_n \colon V_h \rightarrow V_h$ defined by
\begin{align}\label{fafaf}
F_n : = r(k\mathcal{A}_h)^n - \mathcal{S}_h(t_n).
\end{align}
As shown in \cite[Theorem 4.2]{larsson1991finite}, we have
\begin{equation}\label{DPf4}
\|F_n \mathcal{A}_h^{-1} \|_{\mathcal{L}( \dot{H}^0 \times \dot{H}^0 )  } \leq Ck,\; \text{ for }\;  n \geq 1.
\end{equation}
Before proceeding further, we also observe that
\begin{equation}
P_1 F_n \mathbf{P}_h w_0 = P_1 M_hF_n\mathbf{P}_hw_0, \quad \;M_h=\left[\begin{array}{c c}I &0 \\0
&T_h\end{array}\right],
\end{equation}
and thus $\| P_1 F_n \mathbf{P}_h w_0\| = \| P_1 M_hF_n\mathbf{P}_hw_0 \|
\leq \| M_hF_n\mathbf{P}_hw_0 \|_{L_2 (\mathcal{D}) \times L_2 (\mathcal{D})} $.
Accordingly it suffices to focus on the estimate of
$\| M_hF_n\mathbf{P}_hw_0 \|_{L_2 (\mathcal{D}) \times L_2 (\mathcal{D})} $.
By denoting $\widetilde{F}_n := M_h F_n M_h^{-1}$, one can write
\begin{align}\label{DPf6}
M_hF_n\mathbf{P}_hw_0 = M_h F_n M_h^{-1} M_h\mathbf{P}_hw_0 = \widetilde{F}_n M_h \mathbf{P}_hw_0.
\end{align}
Further, using eigenfunction expansions one can easily check that
\begin{equation}
\widetilde{F}_n=r(k\widetilde{\mathcal{A}}_h)^n-e^{-t_n\widetilde{\mathcal{A}}_h},
\; \text{ where }\;
\widetilde{\mathcal{A}}_h=M_h\mathcal{A}_hM_h^{-1}=\left[\begin{array}{c c}0 &-A_h \\I
&\alpha A_h\end{array}\right].
\end{equation}
Since $\widetilde{\mathcal{A}}_h$ shares the same eigenvalues as $\mathcal{\mathcal{A}}_h$, \eqref{DDPf5},\eqref{DPf4} and the stability property $\|e^{-t\widetilde{\mathcal{A}}_h}\|_{\mathcal{L}( \dot{H}^0 \times \dot{H}^0 )}\leq C$ also hold for $\widetilde{\mathcal{A}}_h$.
These facts enable us to derive from \eqref{DPf6} that
\begin{equation}
\begin{split}
\| M_hF_n\mathbf{P}_hw_0 \|_{L_2 (\mathcal{D}) \times L_2 (\mathcal{D})}
& =
\|\widetilde{F}_n\widetilde{\mathcal{A}}_h^{-1}\widetilde{\mathcal{A}}_hM_h\mathbf{P}_hw_0\|_{L_2 (\mathcal{D}) \times L_2 (\mathcal{D})}
\leq
Ck \|\widetilde{\mathcal{A}}_hM_h\mathbf{P}_hw_0\|_{L_2 (\mathcal{D}) \times L_2 (\mathcal{D})}
\\ & =
Ck\bigg\|\left[\begin{array}{c c}0 &-I \\I
&\alpha  I\end{array}\right]\mathbf{P}_hw_0\biggl\|_{L_2 (\mathcal{D}) \times L_2 (\mathcal{D})}
\leq C k ( \| \mathcal{P}_hu_0 \| + \| \mathcal{P}_hv_0 \| ),
\end{split}
\end{equation}
and that
\begin{equation}
\begin{split}
\| M_hF_n\mathbf{P}_hw_0 \|_{L_2 (\mathcal{D}) \times L_2 (\mathcal{D})}
& =
\|\widetilde{F}_nM_h\mathbf{P}_hw_0\|_{L_2 (\mathcal{D}) \times L_2 (\mathcal{D})}
\leq
C\|M_h\mathbf{P}_hw_0\|_{L_2 (\mathcal{D}) \times L_2 (\mathcal{D})}
\\ &
\leq
C( \| \mathcal{P}_hu_0 \| + \| A_h^{-1} \mathcal{P}_hv_0\|).
\end{split}
\end{equation}
Then by interpolation, we obtain
\begin{align}
\| P_1 F_n \mathbf{P}_h w_0\|
\leq
\| M_hF_n\mathbf{P}_hw_0 \|_{L_2 (\mathcal{D}) \times L_2 (\mathcal{D})}
\leq Ck^{\frac s 2}(\|A_h^{\frac s 2}\mathcal{P}_hu_0\|
+
\| A_h^{\frac{ s - 2 }{ 2 } }\mathcal{P}_hv_0\|),\; for\; 0\leq s\leq 2,
\end{align}
which, after assigning $s = 1$ and using \eqref{asr3}, implies \eqref{eq:thm-U-full-estimate1}, that is,
\begin{align}
\| P_1 (r(k\mathcal{A}_h)^n - \mathcal{S}_h(t_n)) \mathbf{P}_h w_0 \|
\leq
C k^{\frac 1 2} \big( \|A_h^{\frac 1 2}\mathcal{P}_hu_0\| + \|A_h^{-\frac 1 2}\mathcal{P}_h v_0 \| \big)
\leq
Ck^{\frac{1}{2}}(|u_0|_1+|v_0|_{-1}).
\end{align}

With regard to (\ref{eq:error-full-v}), the case $q=0$ can be directly obtained by the stability of $e^{-t\mathcal{A}}$ and $r(k\mathcal{A}_h)^n$, and  the case $q=2$ is available in \cite[Theorem 5.1]{larsson1991finite}. Again, the interpolation gives (\ref{eq:error-full-v}).
$\square$

Similarly to the semi-discrete problem as before, we expect a time discrete analogue of
Theorem \ref{thm:error-FEM2}, which requires weaker regularity assumption on $v_0$
than Theorem \ref{thm:error-FEM3} does.  However, this is not an easy job
and, as one can see below, the proof becomes much more involved.
First, we need further regularity results of the linear strongly damped wave equation \eqref{eq:linear-determ-eq}.
\begin{lemma}\label{lem:appen1}
Let $u(t)$ be the solution of the strongly damped wave equation (\ref{eq:linear-determ-eq}), then it holds that
\begin{align}
\int_0^t |u_{tt}(s)|_\beta^2 \, \mathrm{d}s
\leq &
 C(|u_0|_{\beta+2}^2+|v_0|_{\beta+1}^2), \quad  \beta \in \mathbb{R},
    \label{appenn5}
    \\
\int_0^t s^2|u_{ttt}(s)|_{\beta}^2 \, \mathrm{d}s
\leq &
C(|u_0|_{\beta+2}^2+|v_0|_{\beta+1}^2),\quad  \beta \in \mathbb{R},
\label{appen4}
\\
\int_0^t s^2|u_{tt}(s)|_\beta^2 \, \mathrm{d}s
\leq &
 C(|u_0|_{\beta}^2+|v_0|_{\beta-1}^2), \quad  \beta \in \mathbb{R}.
\label{appen5}
\end{align}
\end{lemma}
{\it  Proof of Lemma \ref{lem:appen1}.}
In order to prove (\ref{appenn5}), we multiply both sides of (\ref{eq:linear-determ-eq}) by  $A^\beta u_{tt}$ to obtain
\begin{align}
|u_{tt}(s)|_\beta^2
+
\frac{\alpha}{2} \frac{ \dd}{ \dd s} |u_t(s)|_{\beta+1}^2=-(Au(s),A^\beta u_{tt}(s)).
\end{align}
Integration over $[0, t]$ and the Cauchy-Schwarz inequality yield
\begin{align}
\int_0^t|u_{tt}(s)|_\beta^2\,\dd s
+
\frac{\alpha}{2}|u_t(t)|_{\beta+1}^2
&\leq
\frac{\alpha}{2}|v_0|_{\beta+1}^2
+
\frac{1}{2}\int_0^t|u(s)|_{\beta+2}^2\,\dd s
+
\frac{1}{2}\int_0^t|u_{tt}(s)|_\beta^2\,\dd s.
\end{align}
Further, using  (\ref{eq:S-GROUP-SP-Regu2}) with $\rho=\beta+2, \sigma = 1$ gives
\begin{align}
\int_0^t|u_{tt}(s)|_\beta^2\,\dd s
+
\alpha|u_t(t)|_{\beta+1}^2
&\leq
    \alpha|v_0|_{\beta+1}^2
    +
    \int_0^t|u(s)|_{\beta+2}^2\,\dd s
\leq
 C (|u_0|_{\beta+2}^2+|v_0|_{\beta+1}^2).
\end{align}
To confirm \eqref{appen4},  we differentiate \eqref{eq:linear-determ-eq} with respect to $t$
and multiply both sides by $s^2 A^{\beta} u_{ttt}$ to get
\begin{equation}
s^2 ( u_{ttt}, A^{\beta} u_{ttt} ) + \alpha s^2 ( A u_{tt}, A^{\beta} u_{ttt} ) + s^2( Au_t, A^{\beta} u_{ttt} ) = 0,
\end{equation}
which can be equivalently written as
\begin{equation}
\begin{split}
 s^2(u_{ttt},A^\beta u_{ttt})
  +
   \frac{\alpha}{2} \frac{ \dd}{ \dd s} s^2(A u_{tt},A^\beta u_{tt})
  =
 - \frac{ \dd}{\dd s} s^2(A u_t,A^\beta u_{tt})
 +
 (s (\alpha+s)Au_{tt}+2sAu_t,A^\beta u_{tt}).
\end{split}
\end{equation}
Integration over $[0, t]$ and using the Cauchy-Schwarz inequality suggest that
\begin{equation*}
\begin{split}
\int_0^t s^2 |u_{ttt}(s)|_{\beta}^2\,\dd s
+
\frac{\alpha }{2} t^2 | u_{tt}(t)|_{\beta+1}^2
\leq
   t^2|u_t(t)|_{\beta+1}|u_{tt}(t)|_{\beta+1}
+
   \int_0^t\big[(\alpha +2s)s|u_{tt}(s)|_{\beta+1}^2
+
| u_{t}(s)|_{\beta+1}^2\big] \, \dd s.
\end{split}
\end{equation*}
Applying Gronwall's inequality and taking Lemmas \ref{lem:S-GROUP-spatial-regularity}, \ref{lem:linear-determin-spatial}  into consideration show that
\begin{align}\label{uttt}
\int_0^t s^2|u_{ttt}(s)|_{\beta}^2\,\dd s
+
t^2| u_{tt}(t)|_{\beta+1}^2
\leq
& C \bigg(t^2|u_t(t)|_{\beta+1}|u_{tt}(t)|_{\beta+1}
 +
 \int_0^t \big[ s |u_{tt}(s)|_{\beta+1}^2 +| u_{t}(s)|_{\beta+1}^2 \big]\,\dd s\bigg)
\nonumber\\
\leq &
C(|u_0|_{\beta+2}^2
+
|v_0|_{\beta+1}^2).
\end{align}
This validates (\ref{appen4}).
%
For the estimate of (\ref{appen5}), we, similarly as before, differentiate \eqref{eq:linear-determ-eq} with respect to
$t$ and multiply both sides by $s A^{ \frac{\beta}{2} -1 } $ to acquire
\begin{equation}
 s \alpha A^{ \frac{\beta}{2} } u_{tt}  = - s A^{ \frac{\beta}{2} -1 } u_{ttt} - s A^{ \frac{\beta}{2} } u_{t}.
\end{equation}
Squaring both sides before integration over $[0, t]$ and combining  (\ref{eq:S-GROUP-SP-Regu2})
and (\ref{appen4})  lead us to
\begin{align}
\alpha^2 \int_0^t s^2 |u_{tt}(s)|_\beta^2 \, \dd s
  \leq
      2 \int_0^t s^2 |u_{ttt}(s)|_{\beta-2}^2 \, \dd s
+
2 \int_0^t s^2 |u_{t}(s)|_\beta^2 \, \dd s
\leq
    C(|u_0|_{\beta}^2 + |v_0|^2_{\beta-1}).
\end{align}
This completes the proof of this lemma.
$\square$

%
Now we are ready to formulate the time discrete analogue of Theorem \ref{thm:error-FEM2}.
%
%
\begin{theorem}\label{thm:error-FEM4}
Under the assumptions of Theorem \ref{thm:error-FEM3}, it holds that
\begin{align}
\|P_1F_{kh}^nw_0\| +
\left( k\sum_{m=1}^n\|P_2F_{kh}^mw_0\|^2 \right)^{\frac{1}{2}}
& \leq
C(h^\beta+k^{\frac{\beta}{2}})(|u_0|_\beta+|v_0|_{\beta-1}),\;0\leq \beta\leq 2,
\label{eq:error-full-u-v} \\
\|P_2F_{kh}^nw_0\| & \leq C(h^q+k^{\frac{q}{2}})t_n^{-\frac{q-s}{2}}(|u_0|_{s+2}+|v_0|_s), \;
0\leq s\leq q\leq 2.
\label{eq:error-full-vv}
\end{align}
\end{theorem}
{\it Proof of Theorem \ref{thm:error-FEM4}.}
Since the proof of \eqref{eq:error-full-vv} is easy, we do this first.  The case $s=q$ with $s\in [0,2]$ is a direct consequence of (\ref{eq:error-full-v}).
In addition, the case $s=0, q=2$ can be found in  \cite[Theorem 5.4]{larsson1991finite}.
Similarly to the proof of \eqref{eq:error-FEM-v2},
the desired intermediate case is obvious by interpolation.
In what follows, we focus on the proof of (\ref{eq:error-full-u-v}). Note first that the full-discrete weak variational form of $(\ref{dspr})$ is to find $(U^n,V^n)'\in V_h\times V_h$ such that
\begin{equation} \label{DDPf6}
\begin{split}
(\overline{\partial}_n V^n,\chi_1)+\alpha(\nabla V^n,\nabla\chi_1)+(\nabla U^n,\nabla \chi_1)=0,\;\forall \chi_1\in V_h, \\
(\nabla \overline{\partial}_n U^n,\nabla\chi_2)-(\nabla V^n,\nabla\chi_2)=0,\;\forall \chi_2\in V_h,
\end{split}
\end{equation}
where $\overline{\partial}_n V^n := \frac{V^n-V^{n-1}}{k}$.
Once again, we use interpolation arguments to obtain (\ref{eq:error-full-u-v}).
For the case $\beta=0$, setting $\chi_1=T_hV^n$,
$\chi_2=T_hU^n$ in (\ref{DDPf6}) and adding the resulting two equations give
\begin{equation} \label{eq:U-V-eq}
\left( \overline{\partial}_n U^n, U^n\right)
+
\left( \overline{\partial}_n V^n, T_hV^n\right)
+
\alpha\left( V^n, V^n\right)=0.
\end{equation}
Observing that  $\left( \overline{\partial}_n U^n, U^n\right) \geq \tfrac{1}{2k}  \left( \| U^n \|^2 - \| U^{n-1} \|^2 \right) $
and $\left( \overline{\partial}_n V^n, T_h V^n\right) \geq \tfrac{1}{2k}  \big( \|T_h^{\frac 1 2} V^n\|^2 - \|T_h^{\frac 1 2} V^{n-1}\|^2 \big) $  by the Cauchy-Schwarz inequality, we derive from \eqref{eq:U-V-eq} that
\begin{align}\label{DDPf8}
\|U^n\|^2-\|U^{n-1}\|^2 + \|T_h^{\frac 1 2} V^n\|^2 - \|T_h^{\frac 1 2} V^{n-1}\|^2
+ 2\alpha  k  \| V^n \|^2
\leq 0 .
\end{align}
By summation on $n$ and noting $U^0 = \mathcal{P}_h u_0, V^0 = \mathcal{P}_h v_0$, we deduce that
\begin{equation}\label{DDPf9}
\begin{split}
\|U^n\|^2 + \|T_h^{\frac{1}{2}}V^n\|^2 + 2\alpha\sum_{m=1}^nk\|V^m\|^2 \leq
 \|\mathcal{P}_hu_0\|^2+\|T_h^{\frac{1}{2}}\mathcal{P}_hv_0\|^2
\leq C(\|u_0\|^2+|u_0|_{-1}^2).
\end{split}
\end{equation}
Now we only need to bound $\|u(t_n) \|^2 + \sum_{m=1}^n k \|v(t_m)\|^2$ before we can show \eqref{eq:error-full-u-v} for $\beta=0$.
Observing
\begin{equation} \label{eq:v-expression}
v(t_m)
=
\frac{1}{k} \int_{t_{m-1}}^{t_m}(s-t_{m-1})v_t(s) \,\dd s
 +
\frac{1}{k} \int_{t_{m-1}}^{t_m} v(s) \,\dd s,
\end{equation}
due to integration by parts, we additionally use (\ref{L-Determin-Spa1}) and (\ref{appen5}) to derive
\begin{equation} \label{FDDPf12}
\begin{split}
\sum_{m=1}^n k \|v(t_m)\|^2 & \leq
2 \sum_{m=1}^n k
\bigg(
   \bigg\|
         \frac{1}{k} \int_{t_{m-1}}^{t_m}(s-t_{m-1})v_t(s) \,\dd s
   \bigg\|^2
+
    \bigg\|
    \frac{1}{k} \int_{t_{m-1}}^{t_m} v(s) \,\dd s
    \bigg\|^2
\bigg)
\\
& \leq 2 \int_0^{t_n} s^2 \|v_{t}(s)\|^2 \,\dd s
+
2 \int_0^{t_n} \| v(s) \|^2 \,\dd s
\leq
C( \| u_0 \|^2 + |v_0|^2_{-1}),
\end{split}
\end{equation}
where at the second step we used $0\leq s - t_{m-1} \leq s$.
This together with \eqref{eq:S-GROUP-SP-Regu2} and (\ref{DDPf9}) verifies \eqref{eq:error-full-u-v} for $\beta=0$.
Next we validate \eqref{eq:error-full-u-v} for $\beta=2$. Similarly to \eqref{eq:error-decomposition},
we introduce some notations as
\begin{equation} \label{eq:theta-notations}
\begin{split}
&\theta_1^n=U^n-\mathcal{R}_hu(t_n),\;\rho_1^n=(\mathcal{R}_h-I)u(t_n),
\\
&\theta^n_2=V^n-\mathcal{P}_hv(t_n),\;\rho_2^n=(\mathcal{P}_h-I)v(t_n).
\end{split}
\end{equation}
Combining (\ref{sspdr4}) and (\ref{DDPf6}) yields
\begin{equation}
\begin{split}
\big(\overline{\partial}_n V^n-v_t(t_n),  \chi_1\big)
+
\alpha\big(\nabla \big(V^n-v(t_n) \big), \nabla \chi_1\big)
+\big(\nabla( U^n-u(t_n)),\nabla \chi_1\big)=0, \; \forall \chi_1\in V_h,
\\
\big(\nabla (\overline{\partial}_n U^n-u_t(t_n)), \nabla \chi_2\big)
-
\big(\nabla( V^n-v(t_n)),\nabla \chi_2\big)=0, \;\forall \chi_2 \in V_h.
\end{split}
\end{equation}
Taking the definitions of $\mathcal{P}_h$ and $\mathcal{R}_h$ into account  and plugging the notations
proposed in \eqref{eq:theta-notations} show
\begin{align}
\left( \overline{\partial}_n\theta_2^n, \chi_1 \right)
+
\alpha(\nabla\theta_2^n,\nabla \chi_1)
+
(\nabla\theta_1^n,\nabla \chi_1)
=
( \rho_4^n,  \chi_1)
-
\alpha(\nabla\rho_2^n,\nabla \chi_1),
\; \forall \chi_1\in V_h,
\label{DDPf13} \\
(\nabla\overline{\partial}_n\theta_1^n,\nabla\chi_2)
-
(\nabla \theta_2^n,\nabla \chi_2)
=
(\nabla \rho_3^n , \nabla \chi_2)
+
(\nabla \rho_2^n,\nabla \chi_2),
\; \forall \chi_2\in V_h,
\label{DDPf14}
\end{align}
where further notations were also introduced:
\begin{equation}
\rho_3^n := u_t(t_n)-\overline{\partial}_nu(t_n) \quad \text{ and } \quad
\rho_4^n := v_t(t_n)-\overline{\partial}_nv(t_n).
\end{equation}
As in the proof of (\ref{eq:semi-discrete-error1}), setting $\chi_1=T_h\theta_2^n $, $\chi_2=T_h\theta_1^n$ in (\ref{DDPf13})-(\ref{DDPf14}) and adding together give
\begin{align}
 \left(\overline{\partial}_n\theta_2^n, T_h \theta_2^n \right)
+
( \overline{\partial}_n\theta_1^n,\theta_1^n)
+\alpha \| \theta_2^n \|^2
 & =
(T_h\mathcal{P}_h \rho_4^n, \theta_2^n)
-
\alpha(\mathcal{R}_h\rho_2^n,\theta_2^n)
+
(\mathcal{R}_h \rho_3^n ,  \theta_1^n)
+
(\mathcal{R}_h\rho_2^n,\theta_1^n)
\nonumber\\
& =
(\tfrac{1}{ \sqrt{\alpha} } T_h\mathcal{P}_h \rho_4^n - \sqrt{\alpha} \mathcal{R}_h\rho_2^n, \sqrt{\alpha} \theta_2^n)
+
(\mathcal{R}_h \rho_3^n + \mathcal{R}_h\rho_2^n, \theta_1^n ).
\end{align}
Using the facts $ \left(  \overline{\partial}_n\theta_1^n, \theta_1^n \right)
\geq \frac{1}{2k}(\|\theta_1^n\|^2- \|\theta_1^{n-1}\|^2) $ and
$ \left(  \overline{\partial}_n\theta_2^n, T_h \theta_2^n \right)
\geq \frac{1}{2k}(\| T_h^{\frac12} \theta_2^n\|^2 - \| T_h^{\frac12}  \theta_2^{n-1}\|^2) $ shows
\begin{align}
& \tfrac{1}{2k}\Big( \|\theta_1^n\|^2
-
\|\theta_1^{n-1}\|^2
+
\|T_h^{\frac12 }\theta_2^n\|^2
-
\|T_h^{\frac12}\theta_2^{n-1}\|^2 \Big)
+
\alpha \|\theta_2^n\|^2
\nonumber\\
& \quad \leq
\tfrac{1}{\alpha}\|T_h\mathcal{P}_h \rho_4^n\|^2
+
(\alpha + 1 )\|\mathcal{R}_h\rho_2^n\|^2
+
\tfrac \alpha 2 \|\theta_2^n\|^2
+
\|\mathcal{R}_h\rho_3^n\|^2
+
\tfrac12 \|\theta_1^n\|^2.
\end{align}
Hence, by summation and detecting that $\theta_2^0=V^0-\mathcal{P}_hv_0=0$ we infer
\begin{equation*}
\begin{split}
\|T_h^{\frac{1}{2}}\theta_2^n\|^2
+\alpha k\sum_{m=1}^n
\|\theta_2^m\|^2+\|\theta_1^n\|^2
\leq&
Ck\sum_{m=1}^n
\left(
\|\mathcal{R}_h\rho_3^m\|^2
+\left\|T_h\mathcal{P}_h\rho_4^m\right\|^2
+\left\|\mathcal{R}_h\rho_2^m\right\|^2
+ \|\theta_1^m\|^2
\right)
+\| \theta_1^0\|^2.
\end{split}
\end{equation*}
%
%
Applying the discrete Gronwall inequality helps us to get
\begin{align}\label{DDPf15}
\|T_h^{\frac{1}{2}}\theta_2^n\|^2
+
\alpha k\sum_{m=1}^n\|\theta_2^m\|^2
+
\|\theta_1^n\|^2
\leq&
 Ck\sum_{m=1}^n \! \left(\|\mathcal{R}_h\rho_3^m\|^2
+
\left\|T_h\mathcal{P}_h\rho_4^m\right\|^2
+
\left\|\mathcal{R}_h\rho_2^m\right\|^2\right) \!
+
Ch^4| u_0|_2^2,
\end{align}
where $ \| \theta_1^0 \| = \| U_0 - \mathcal{R}_h u_0 \| = \| \mathcal{P}_h(I-\mathcal{R}_h)u_0 \| \leq C h^2 |u_0|_2 $. In the sequel  we will estimate the remaining three terms separately. Note first that $\rho_3^m$ admits the following expression
\begin{align}\label{DDPf16}
\rho_3^m = u_t(t_m)-\overline{\partial}_mu(t_m)=\frac{1}{k}
 \int_{t_{m-1}}^{t_m}(s-t_{m-1})u_{tt}(s) \, \text{d} s.
\end{align}
This together with (\ref{appenn5}) and (\ref{appen5}) guarantees
\begin{align}
\label{DDPf17}
k \sum_{m=1}^n \|& \mathcal{R}_h\rho_3^m \|^2
\leq
2 k \sum_{m=1}^n \| (I-\mathcal{R}_h)\rho_3^m\|^2
   + 2k \sum_{m=1}^n \| \rho_3^m \|^2
\nonumber \\
& \leq
C k\sum_{m=1}^n \bigg( h^4 \bigg| \frac{1}{k} \int_{t_{m-1}}^{t_m} ( s - t_{m-1} ) u_{tt} ( s ) \, \dd s \bigg|_2^2
+
 \bigg\| \frac{1}{k} \int_{t_{m-1}}^{t_m} (s-t_{m-1})u_{tt}(s) \, \dd s \bigg\|^2 \bigg) \nonumber\\
&\leq
     C(h^4+k^2) \int_0^{t_n}( s^2 | u_{tt}(s) |_2^2 +\|u_{tt}(s)\|^2)\,\dd  s
     \leq
       C(h^4+k^2)(|u_0|_2^2+|v_0|_1^2),
\end{align}
where the fact $s - t_{m - 1} \leq s$ was used.
Likewise,  noting that $T_h\mathcal{P}_h=T_h$ and that $\rho_4^m$ has the same expression as \eqref{DDPf16} with $u$ replaced by $v$ yields
\begin{align} \label{DDPf18}
k \sum_{m=1}^n \| T_h\mathcal{P}_h\rho_4^m \|^2
\leq&
2k \sum_{m=1}^n \|(T-T_h) \rho_4^m\|^2
+
2k \sum_{m=1}^n \|T\rho_4^m\|^2
\nonumber \\
\leq &
    C k\sum_{m=1}^n \bigg( h^4 \bigg\| \frac{1}{k}\int_{t_{n-1}}^{t_n}(s-t_{n-1})v_{tt}(s) \,\dd s\bigg \|^2
+
  \bigg \|  \frac{1}{k}  \int_{t_{n-1}}^{t_n}   (s-t_{n-1})Tv_{tt}(s) \,\dd s  \bigg\|^2 \bigg)  \nonumber\\
\leq&
   Ch^4\int_0^{t_n}s^2\|u_{ttt}(s)\|^2\, \dd s + Ck^2 \int_0^{t_n} (\|u_{tt}(s)\|^2+\|u_t(s)\|^2) \,\dd s\nonumber\\
\leq&
    C(h^4+k^2)(|u_0|_2^2+|v_0|_1^2),
\end{align}
where we also used \eqref{L-Determin-Spa1}, (\ref{dadlal}), \eqref{appenn5}, \eqref{appen4}  and
$Tu_{ttt} = - \alpha u_{tt} - u_t$.
Using similar arguments as before and taking \eqref{IIII1}, \eqref{L-Determin-Spa1}, \eqref{appen5} and \eqref{eq:v-expression} into account one can show that
\begin{align} \label{DDPf19}
   k\sum_{m=1}^n \|\mathcal{R}_h\rho_2^m\|^2
= &
   k\sum_{m=1}^n\|  \mathcal{P}_h (\mathcal{R}_h-I)  v(t_m)\|^2
\leq
   Ch^4k\sum_{m=1}^n|v(t_m)|^2_2
   \nonumber\\
\leq&
     Ch^4k\sum_{m=1}^n  \biggl|  \frac{1}{k}  \int_{t_{m-1}}^{t_m} (s-t_{m-1})v_{t}(s) \,\dd s
+
   \frac{1}{k}  \int_{t_{m-1}}^{t_m}  v(s)\,\dd s \biggl|_2^2
\nonumber\\
\leq&
     Ch^4\left( \int_0^{t_n} s^2 | u_{tt}(s) |_2^2 \, \dd s
+
\int_0^{t_n} | u_t(s) |_2^2 \, \dd s \right)
\leq
    Ch^4( |u_0|_2^2 + |v_0|_1^2 ).
\end{align}
Analogously, one can achieve
\begin{equation} \label{eq:rho2-sum}
k \sum_{m=1}^n \| \rho_2^m\|^2 + \|\rho_1^n\|^2  \leq C(h^4+k^2)(|u_0|_2^2+|v_0|_1^2).
\end{equation}
Finally,  plugging (\ref{DDPf17})-(\ref{DDPf19})  into (\ref{DDPf15}) and considering \eqref{eq:rho2-sum}  help us to get
\begin{align}
k\sum_{m=1}^n \| V^m - v(t_m) \|^2 + \| U^n - u(t_n) \|^2 \leq C( h^4 + k^2 ) ( |u_0|_2^2 + |v_0|_1^2 ).
\end{align}
The intermediate cases follow by interpolation.
%
$\square$
\section{Finite element method for the stochastic problem}
\label{sec:stochastic-main-result}
This section is devoted to the finite element approximation of the stochastic problem (\ref{stojab}).
The convergence analysis relies on regularity properties of the mild solution of (\ref{stojab})
as well as error estimates obtained in section \ref{sec:Error-estimate-determin}.
 %
%

\subsection{Spatial semi-discretization}
In this subsection, we shall follow notations introduced in  section \ref{sec:Error-estimate-determin} and
analyze the semidiscrete finite element approximation of (\ref{stojab}).
Let $V_h$ be the finite element  space defined in the previous section.
The semidiscrete approximation of (\ref{stojab}) is to find $X_h(t)=(u_h(t), u_{h,t}(t))'\in V_h\times V_h$ such that
\begin{align} \label{sttodae2}
 \textmd{d}X_h(t)
+
\mathcal{A}_hX_h(t)\textmd{d}t
=
\mathbf{P}_h\mathbf{F}(X_h(t)) \, \textmd{d}t
+
\mathbf{P}_h\mathbf{B}
\, \textmd{d}W(t),
\;
  \text{ in }\; t \in (0, T],\;X_h(0)
=
  \mathbf{P}_hX_0,
\end{align}
or in the mild form
\begin{align}\label{eq:mild-FEM}
X_h(t)=\mathcal{S}(t)\mathbf{P}_hX_0+\int_0^t\mathcal{S}_h(t-s)\mathbf{P}_h\mathbf{F}(X_h(s)) \, \textmd{d} s
+
\int_0^t\mathcal{S}_h(t-s)\mathbf{P}_h\mathbf{B}\, \textmd{d} W(s), \; t\in[0,T].
\end{align}

The first main convergence result is as follows.
\begin{theorem}\label{thm:conteh1}
Let Assumptions \ref{ass:F}-\ref{ass:initial-value} hold with $\gamma \in [0, 1]$ and let the setting in the beginning of section \ref{sec:Error-estimate-determin} be fulfilled. Let $(u(t), u_t(t))'$ and
$(u_h(t), u_{h,t}(t))'$ be  the mild solutions of the problems (\ref{stojab}) and (\ref{sttodae2}), respectively.
Then for all $t \in [0, T]$ it holds that
\begin{equation}
\|u(t)-u_h(t)\|_{L^2(\Omega;\dot{H}^0)}  \leq C h^{1+\gamma}
\big(
     1 + \|\varphi\|_{L^2(\Omega;\dot{H}^{\gamma+1})} + \|\psi\|_{L^2(\Omega; \dot{H}^{\gamma -1} ) }
\big). \label{thuh1}
\end{equation}
If additionally $\psi  \in L^2( \Omega; \dot{H}^{\gamma} )$, then for all $t \in [0, T]$ it holds that
\begin{equation}
\|u_t(t)-u_{h,t}(t)\|_{L^2(\Omega;\dot{H}^0)}  \leq C h^{\gamma}
\big(
1 + \|\varphi\|_{L^2(\Omega; \dot{H}^{ \gamma })} + \|\psi\|_{L^2(\Omega; \dot{H}^{\gamma})}
\big).
\label{thuh2}
\end{equation}
\end{theorem}
\begin{remark}
When Assumption \ref{ass:initial-value} is fulfilled with $\gamma\in[-1,0)$, the problem
(\ref{eq:intro-SSDVE}) can admit a mild solution $\{ u(t) \}_{t \in [0, T]}$ that exhibits a
positive order of regularity. However, in this situation the Wiener process
takes values in $\dot{H}^{\delta}$ for $\delta = \gamma - 1 < -1$, which destroys the well-posedness of $\mathbf{P}_h$,  as also explained in Remark \ref{remark:banqdr2}.
Therefore, throughout this section we restrict ourselves to $\gamma \in [0, 1]$.
\end{remark}
{\it Proof of Theorem \ref{thm:conteh1}.}
Subtracting \eqref{eq:mild-SPDE} from \eqref{eq:mild-FEM} gives
\begin{align}\label{diffener}
X_h(t)&-X(t)
=
\big(  \mathcal{S}_h(t)\mathbf{P}_h-\mathcal{S}(t)  \big)X_0
+
  \int_0^t \big(   \mathcal{S}_h(t-s)\mathbf{P}_h  -\mathcal{S}(t-s)  \big)   \mathbf{F}(X(s)) \,\dd s
   \nonumber\\
&+
\int_0^t \mathcal{S}_h(t-s)\mathbf{P}_h  \Big( \mathbf{F}\big(X_h(s)  \big) -\mathbf{F}\big(X(s)  \big)   \Big)\,\dd s
+
\int_0^t \big(  \mathcal{S}_h(t-s)\mathbf{P}_h -\mathcal{S}(t-s) \big) \mathbf{B}\,\dd W(s)
\nonumber\\
: =&J_1+J_2+J_3+J_4.
\end{align}
Recalling $u(t)-u_h(t)=P_1(X(t)-X_h(t))$, we require to bound $P_1 J_i, i =1,2,3,4$.
For the term $P_1J_1$,  a combination with (\ref{eq:error-FEM-u}) and (\ref{eq:error-FEM-u-v}) enables us to claim that, for $\beta\in [0, \gamma]$, $i\in\{0,1\}$,
\begin{align}\label{eq:unified-results-F(t)}
\|P_1F_h(t)w_0\|\leq Ch^{\beta+i} (|u_0|_{\beta+i}+|v_0|_{\beta-1}),
\end{align}
which together with  Assumption \ref{ass:initial-value} leads to
  \begin{align}\label{uuh1}
\|P_1J_1\|_{L^2(\Omega;\dot{H}^0)}
= &
\|P_1 F_h (t)X_0\|_{L^2(\Omega;\dot{H}^0)}
\leq Ch^{\beta+i}
 \big(
     \|\varphi\|_{L^2(\Omega;\dot{H}^{\beta+i})}
     +
     \|\psi\|_{L^2(\Omega; \dot{H}^{\beta -1} )}
\big).
\end{align}
Similarly, using \eqref{eq:error-FEM-u} and \eqref{eq:F-moment-bound} shows
\begin{align}
\|P_1J_2\|_{L^2(\Omega;\dot{H}^0)}
& \leq
\int_0^t \|P_1 F_h (t -s) \mathbf{F}(X(s))\|_{L^2(\Omega;\dot{H}^0)} \,\dd s
\leq
C h^2 \int_0^t \| F(u(s)) \|_{L^2(\Omega;\dot{H}^0)} \,\dd s
\nonumber \\
&
\leq
   C h^2  \big(1 + \|\varphi\|_{ L^2(\Omega;\dot{H}^0 ) } + \|\psi\|_{L^2(\Omega; \dot{H}^{-2})}  \big).
\end{align}
To bound $P_1J_3$, we combine the stability of $\mathcal{S}_h(t) \mathbf{P}_h$ in $ \dot{H}^0 \times \dot{H}^0$ with Assumption \ref{ass:QWiener} to derive
\begin{align}
\|P_1J_3\|_{L^2(\Omega;\dot{H}^0)}
 &\leq
  \int_0^t \left\|
   P_1 \mathcal{S}_h(t-s)\mathbf{P}_h
       \Big( \mathbf{F}\big(X_h(s)  \big)
           -\mathbf{F}\big(X(s)  \big)   \Big)   \right\|_{L^2(\Omega;\dot{H}^0)}\,\dd s
 \nonumber\\
&\leq
   C\int_0^t   \|  F(u_h(s))-F(u(s))  \|_{L^2(\Omega;\dot{H}^{0})}  \,\dd s
  \leq
      C\int_0^t  \|u_h(s)-u(s)\|_{L^2(\Omega;\dot{H}^0)}  \,\dd s.
\end{align}
Again, using (\ref{eq:error-FEM-u}) and the It\^{o} isometry yields
\begin{equation}\label{uuh1}
\begin{split}
\|P_1J_4\|_{L^2(\Omega;\dot{H}^0)}
=
\left\|\int_0^t P_1 F_h ( t - s ) \mathbf{B} \, \text{d} W(s)\right\|_{L^2(\Omega;\dot{H}^0)}
& =
\left(\int_0^t
    \|P_1 F_h ( t - s ) \mathbf{B}Q^{\frac{1}{2}}\|_{\mathrm{HS}}^2 \, \text{d} s
       \right)^{\frac{1}{2}}
\\ & \leq
    C\sqrt{T} h^{\gamma+1}   \|A^{\frac{\gamma-1}{2}}Q^{\frac{1}{2}}\|_{\mathrm{HS}}.
\end{split}
\end{equation}
Finally,  putting the above estimates together and employing Gronwall's inequality give
\begin{align}\label{eq:error-u-uh}
\|u(t)-u_h(t)\|_{L^2(\Omega;\dot{H}^0)}\leq Ch^{\beta+i}(1+ \|\varphi\|_{L^2(\Omega;\dot{H}^{\beta+i})}
     +
     \|\psi\|_{L^2(\Omega; \dot{H}^{\beta -1} )}), \quad \beta\in[0,\gamma], \, i \in \{0,1\}.
\end{align}
Letting  $\beta=\gamma$, $i=1$ in (\ref{eq:error-u-uh}) hence yields (\ref{thuh1}).
Next, we are to verify (\ref{thuh2}). Following the same notations as before, we need to estimate
$P_2J_i, i = 1, 2, 3, 4$.
Using   (\ref{eq:error-FEM-v}) with $q = \gamma$ gives
  \begin{align}
\|P_2J_1\|_{L^2(\Omega;\dot{H}^0)}
= &
    \|P_2 F_h (t) X_0 \|_{L^2(\Omega;\dot{H}^0)}
\leq Ch^{\gamma}
    \big(
         \|\varphi\|_{L^2(\Omega;\dot{H}^{\gamma})}
+
           \|\psi\|_{L^2(\Omega; \dot{H}^{\gamma})}
\big).
\end{align}
To deal with the term $P_2J_2$, we employ  (\ref{eq:error-FEM-v2}) with $q = \gamma, s = 0$ and
\eqref{eq:F-moment-bound} to arrive at
\begin{align}
\|P_2J_2\|_{L^2(\Omega;\dot{H}^0)}
\leq&
\int_0^t
     \big\| P_2 F_h ( t -s ) \mathbf{F}(X(s))
          \big \|_{L^2(\Omega;\dot{H}^0)}\,\dd s
\leq
    Ch^{\gamma}\int_0^t    (t-s)^{-\frac{\gamma}{2}} \|F(u(s))\|_{L^2(\Omega;\dot{H}^0)}\,\dd s
\nonumber \\
\leq &
    C h^{\gamma} \big(1 + \|\varphi\|_{L^2(\Omega; \dot{H}^{ 0 })}
+
\|\psi\|_{ L^2(\Omega; \dot{H}^{ -2 } ) } \big).
\end{align}
The stability of $ \mathcal{S}_h(t)\mathbf{P}_h $ in $\dot{H}^0 \times \dot{H}^0$, (\ref{eq:error-u-uh}) with $\beta=\gamma$, $i=0$ and Assumption \ref{ass:F}  ensure
\begin{align}
\|P_2J_3\|_{L^2(\Omega;\dot{H}^0)}
&\leq
     \int_0^t \Big\| P_2 \mathcal{S}_h ( t - s ) \mathbf{P}_h
              \Big( \mathbf{F}(X(s)) - \mathbf{F} ( X_h(s) ) \Big) \Big\|_{L^2(\Omega;\dot{H}^0)} \,\dd s
\nonumber \\
&\leq
    C\int_0^t  \|u(s)-u_h(s)\|_{L^2(\Omega;\dot{H}^0)}  \,\dd s
\leq C h^{ \gamma }
\big(
1 + \|\varphi\|_{L^2(\Omega; \dot{H}^{ \gamma })} + \|\psi\|_{L^2(\Omega; \dot{H}^{\gamma - 1})}
\big).
\end{align}
At last, It\^{o}'s isometry and (\ref{eq:error-FEM-u-v}) with $\beta = \gamma$ help us to estimate $P_2J_4$ as follows,
\begin{align}\label{uuh1}
\|P_2J_4\|_{L^2(\Omega;\dot{H}^0)}
=
    \left(  \int_0^t \|P_2 F_h (t -s)  \mathbf{B}Q^{\frac{1}{2}}\|_{\mathrm{HS}}^2 \,\dd s   \right)^{\frac{1}{2}}
 \leq
    Ch^{\gamma}\|A^{\frac{\gamma-1}{2}}Q^{\frac{1}{2}}\|_{\mathrm{HS}}.
\end{align}
Now gathering the estimates of $P_2J_i, i=1,2,3,4$ together gives the estimate of $u_t(t)-u_{h,t}(t)$.
$\square$

\subsection{Full-discretization}
%
Below, we proceed to treat the full-discrete scheme for  \eqref{stojab}.
Let $k$ be the time step-size and  write $t_n=nk$, for $n\geq 1$. We
discretize (\ref{sttodae2}) in time with a linear implicit Euler scheme and the resulting full-discretization
is thus to find $\mathcal{F}_{t_n}$-adapted random variables $X^n=(U^n, V^n)'\in V_h\times V_h$ such that
\begin{align} \label{stofdae2}
 X^n-X^{n-1}+k\mathcal{A}_hX^n=k \mathbf{P}_h\mathbf{F}(X^{n-1})+\mathbf{P}_h\mathbf{B}\Delta W_n,\;
X^0=\mathbf{P}_hX_0,\;in\;\mathcal{D},
\end{align}
where $\Delta W_n : =W(t_n)-W(t_{n-1})$  is the Wiener increment.
%
Now we  state our last convergence result.
\begin{theorem}\label{thm:conteh2}
Let $(u(t), u_t(t))'$ and $(U^n, V^n)'$ be the solutions of \eqref{stojab} and (\ref{stofdae2}), respectively.
If Assumptions \ref{ass:F}--\ref{ass:initial-value} hold with $ \gamma \in [0, 1]$ and the setting in the beginning of section \ref{sec:Error-estimate-determin} holds, then
\begin{align}\label{quansq1}
\|u(t_n)-U^n\|_{L^2(\Omega;\dot{H}^0)} \leq C (h^{\gamma+1}+k^{\frac{ \gamma+1 } {2} } )
\big(
1+\|\varphi\|_{L^2(\Omega;\dot{H}^{\gamma+1})}+\|\psi\|_{L^2(\Omega; \dot{H}^{\gamma - 1} ) }
\big).
\end{align}
If additionally $\psi  \in L^2( \Omega; \dot{H}^{\gamma} )$, then
\begin{equation}
\|u_t(t_n)-V^n\|_{L^2(\Omega;\dot{H}^0)} \leq C (h^{\gamma}+k^{\frac{\gamma}{2}})
\big(
1+\|\varphi\|_{L^2(\Omega;\dot{H}^{\gamma})} + \|\psi\|_{L^2(\Omega; \dot{H}^{\gamma})}
\big).
\label{qudddansq1}
\end{equation}
\end{theorem}
We begin by introducing a crucial ingredient in the following convergence analysis.
\begin{lemma}\label{lemmad}
Suppose that  $ w_0 = (u_0, v_0)' \in \dot{H}^\mu \times \dot{H}^{\mu-1} $ for some $\mu\in [0,2]$. Then
\begin{align}
\sum_{j=0}^{n-1}\int_{t_j}^{t_{j+1}}\|P_2(\mathcal{S}(t_n-s)-\mathcal{S}(t_n-t_j))w_0\|^2
\,\mathrm{d}s
   \leq
      Ck^{\mu}(|u_0|^2_\mu+|v_0|_{\mu-1}^2 ).
\label{shjian1}
\end{align}
\end{lemma}
{\it Proof of Lemma \ref{lemmad}.}
Keep in mind that
\begin{equation}
P_2(\mathcal{S}(t_n-s)-\mathcal{S}(t_n-t_j))w_0 = u_t(t_n-s)-u_t(t_n-t_{j}),
\end{equation}
where by abuse of notation we view $ (u(t),u_t(t))'$ as the solution of the equation (\ref{eq:linear-determ-eq}).
By interpolation, we only need to verify \eqref{shjian1} for the cases $\mu=0$ and $\mu=2$. Using  (\ref{appenn5}) shows
\begin{align}\label{dddalal}
\sum_{j=0}^{n-1}  \int_{t_j}^{t_{j+1}}  \|u_t(t_n-s)-u_t(t_n-t_{j})\|^2\,\dd s
=&
  \sum_{j=0}^{n-1}   \int_{t_{j}}^{t_{j+1}}    \Big \| \int_{t_j}^{s}u_{tt}(t_n-r) \,\dd r  \Big \|^2  \,\dd s
     \nonumber\\
\leq &
   C k^2\int_0^{t_n}\|u_{tt}(t_n-s)\|^2\,\dd s
   \leq
   C k^2(|u_0|_2^2+|v_0|_1^2).
\end{align}
Also, employing (\ref{L-Determin-Spa1}) with $\beta = 0$ and (\ref{FDDPf12}) shows
\begin{align*}\label{dddalaldd}
\begin{split}
\sum_{j=0}^{n-1}\int_{t_{j}}^{t_{j+1}} \|u_t(t_n-s)-u_t(t_n-t_{j})\|^2 \,\dd s
\leq  2k\sum_{j=1}^{n}\|u_t(t_j)\|^2
+
  2\int_0^{t_n} \! \|u_t(t_n-s)\|^2  \,\dd s
   \leq
      C(\|u_0\|^2+|v_0|_{-1}^2).
\end{split}
\end{align*}
This and (\ref{dddalal}) together concludes the proof of this lemma.
$\square$

{\it Proof of Theorem \ref{thm:conteh2}.}
Equivalently, (\ref{stofdae2}) can be reformulated as
\begin{align}\
X^n = r(k\mathcal{A}_h)^n\mathbf{P}_hX_0
+
\sum_{j=0}^{n-1} \int_{t_j}^{t_{j+1}}r(k\mathcal{A}_h)^{n-j}\mathbf{P}_h\mathbf{F}(X^{j})\,\dd s
+
\sum_{j=0}^{n-1}\int_{t_{j}}^{t_{j+1}}r(k\mathcal{A}_h)^{n-j}\mathbf{P}_h \mathbf{B} \,\dd W(s).
\end{align}
Therefore, the difference between $X^n$ and $X(t_n)$ can be decomposed as follows:
\begin{equation}
\begin{split}
X^n -X(t_n) = & (r(k\mathcal{A}_h)^n\mathbf{P}_h-\mathcal{S}(t_n))X_0
+
\sum_{j=0}^{n-1} \int_{t_j}^{t_{j+1}}  r(k\mathcal{A}_h)^{n-j} \mathbf{P}_h
          \left(  \mathbf{F}(X^{j})-\mathbf{F}(X(s))  \right)\,\dd s
      \\
&+
\sum_{j=0}^{n-1}   \int_{t_{j}}^{t_{j+1}}   \left(  r(k\mathcal{A}_h)^{n-j}\mathbf{P}_h-\mathcal{S}(t_{n-j})
                                               \right) \mathbf{F}(X(s)) \,\dd s
                           \\
& +
\sum_{j=0}^{n-1}\int_{t_{j}}^{t_{j+1}}
          \left(  \mathcal{S}(t_{n-j})-\mathcal{S}(t_n-s)  \right) \mathbf{F}(X(s))\,\dd s
          \\
& +
\sum_{j=0}^{n-1}\int_{t_{j}}^{t_{j+1}}  \left(  r(k\mathcal{A}_h)^{n-j}\mathbf{P}_h-\mathcal{S}(t_n-s)
                                          \right) \mathbf{B}\,\dd W(s)
:=
\mathbb{J}_1 + \mathbb{J}_2 + \mathbb{J}_3 + \mathbb{J}_4 + \mathbb{J}_5.
\end{split}
\end{equation}
Since  $U^n-u(t_n)= P_1 ( X^n -X(t_n) ) $,  the estimate $\|u(t_n)-U^n\|_{L^2(\Omega;\dot{H}^0)}$ can be achieved via estimates of $\| P_1\mathbb{J}_j  \|_{L^2(\Omega;\dot{H}^0)}$, $j =1,2,...,5$.
As in (\ref{eq:unified-results-F(t)}), combining (\ref{eq:error-full-u}) and (\ref{eq:error-full-u-v}) implies, for all $\beta\in[0,\gamma]$, $i\in\{0,1\}$,
\begin{align}
\|P_1F_{kh}^nw_0\|\leq C(h^{\beta+i}+k^{\frac{\beta+i}{2}})(|u_0|_{\beta+i}+|v_0|_{\beta-1}).
\end{align}
This immediately leads us to the estimate of $P_1 \mathbb{J}_1$ as follows,
\begin{align}\label{l111es}
\| P_1 \mathbb{J}_1 \|_{L^2(\Omega;\dot{H}^0)}\leq C(h^{\beta+i}+k^{\frac{\beta+i}{2}})(\|\varphi\|_{L^2(\Omega;\dot{H}^{\beta+i})}+\|\psi\|_{L^2(\Omega; \dot{H}^{\beta-1})}),  \quad \beta\in[0,\gamma], \, i \in\{0,1\} .
\end{align}
For $P_1  \mathbb{J}_2$,  we use the stability property \eqref{DDPf5}, the regularity  (\ref{eq:u(t)-u(s).final}) with $\varrho=\beta+i-1$, $i\in\{0,1\}$, and Assumption \ref{ass:F} to get
\begin{align}\label{aldlalg}
\| P_1  \mathbb{J}_2 \|_{L^2(\Omega;\dot{H}^0)}
\leq &
\sum_{j=0}^{n-1} \int_{t_j}^{t_{j+1}} \big\| P_1r(k\mathcal{A}_h)^{n-j}\mathbf{P}_h
           \big( \mathbf{F}(X^{j})-\mathbf{F}(X(s)) \big) \big \|_{L^2(\Omega;\dot{H}^0)} \, \dd s
\nonumber\\
\leq &
     C\sum_{j=0}^{n-1} \int_{t_j}^{t_{j+1}}\|F(U^j)-F(u(s))\|_{L^2(\Omega;\dot{H}^{0})}  \, \dd s
\nonumber\\
\leq &
   C\sum_{j=0}^{n-1} \int_{t_j}^{t_{j+1}}\Big(\|U^j-u(t_j)\|_{L^2(\Omega;\dot{H}^0)}
   +
   \|u(t_j)-u(s)\|_{L^2(\Omega;\dot{H}^0)}\Big) \, \dd s
\nonumber\\
\leq&
 C k \sum_{j=0}^{n-1}\|U^j-u(t_j)\|_{L^2(\Omega;\dot{H}^0)}
+
C k^{\frac{\beta+i}{2}} \big( 1+\|\varphi\|_{L^2(\Omega;\dot{H}^{\beta+i})}
 +
 \|\psi\|_{L^2(\Omega; \dot{H}^{\beta+i-2})} \big).
\end{align}
To bound the term $P_1  \mathbb{J}_3$, we recall (\ref{eq:error-full-u}) with $ \beta = 2 $ and derive that
\begin{align}\label{leq31}
\| P_1  \mathbb{J}_3 \|_{L^2(\Omega;\dot{H}^0)}
\leq  &
 \sum_{j=0}^{n-1}\int_{t_{j}}^{t_{j+1}}
       \| P_1 F_{kh}^{n-j} \mathbf{F}(X(s))  \|_{L^2(\Omega;\dot{H}^0)}  \,\dd s
\leq
      C( h^{2} + k )T \sup_{ s\in [0, T]} \| F(u(s)) \|_{L^2(\Omega;\dot{H}^0)}
\nonumber \\
\leq &
    C( h^{2} + k )
\big( 1+\|\varphi\|_{L^2(\Omega;\dot{H}^{0})}
+
   \|\psi\|_{L^2(\Omega; \dot{H}^{-2})} \big).
\end{align}
In the same spirit as before but employing \eqref{eq:linear-deter-u-time-regularity-S} with $ \mu = 2$ instead we obtain
\begin{equation}
\begin{split}
\| P_1  \mathbb{J}_4 \|_{L^2(\Omega;\dot{H}^0)} &
\leq
\sum_{j=0}^{n-1}\int_{t_{j}}^{t_{j+1}}
    \|P_1(\mathcal{S}(t_{n-j})-\mathcal{S}(t_n-s)) \mathbf{F}(X(s))\|_{L^2(\Omega;\dot{H}^0)} \, \dd s
\\
&\leq
Ck  \big( 1+\|\varphi\|_{L^2(\Omega;\dot{H}^{0} )} + \|\psi\|_{L^2(\Omega; \dot{H}^{-2} ) } \big).
\end{split}
\end{equation}
Now it remains to estimate $ P_1  \mathbb{J}_5 $. Employing (\ref{eq:linear-deter-u-time-regularity-S}), \eqref{eq:error-full-u}  and  It\^{o}'s isometry together promises
\begin{align}\label{adddadd}
\|  P_1  \mathbb{J}_5 \|^2_{L^2(\Omega;\dot{H}^0)}
= &
\sum_{j=0}^{n-1}   \int_{t_{j}}^{t_{j+1}}   \left\|
               P_1\left(r(k\mathcal{A}_h)^{n-j}\mathbf{P}_h-\mathcal{S}(t_n-s)\right) \mathbf{B}Q^{\frac{1}{2}}
                                           \right\|^2_{\mathrm{HS}} \,\dd s
\nonumber\\
\leq &
 2 \sum_{j=0}^{n-1}
     \int_{t_{j}}^{t_{j+1}}
                \left(\| P_1 F_{ kh }^{ n-j } \mathbf{B}Q^{\frac 1 2} \|_{ \mathrm{HS} }^2
+
                   \| P_1(\mathcal{S}(t_{n-j}) - \mathcal{S}(t_n-s))
                                \mathbf{B}Q^{\frac 1 2}  \|_{\mathrm{HS}}^2\right) \,\dd s
\nonumber\\
\leq&
 C(h^{2(\gamma+1)}+k^{\gamma+1})
                    \|A^{\frac{\gamma-1}{2}}Q^{\frac 1 2}\|_{\mathrm{HS}}^2.
\end{align}
Putting the above five estimates together and applying Gronwall's inequality imply that,
\begin{align}\label{eq:error-u-U-general}
\|u(t_n)-U^n\|\leq C(h^{\beta+i}+k^{\frac{\beta+i}{2}}) (1+\|\varphi\|_{L^2(\Omega;\dot{H}^{\beta+i})}+\|\psi\|_{L^2(\Omega; \dot{H}^{\beta-1})}) \:\:  \text{ for }  \beta\in[0,\gamma], \, i \in \{0,1\},
\end{align}
which validates \eqref{quansq1} by taking $\beta = \gamma, i = 1$.
In the sequel we turn our attention to the estimate of $V^n-v(t_n)$.
Using  (\ref{eq:error-full-v})  with $q=\gamma$ suggests
\begin{align}\label{l12es}
\| P_2 \mathbb{J}_1 \|_{L^2(\Omega;\dot{H}^0)}
    \leq C \| P_2 F_{kh}^n X_0 \|_{L^2(\Omega;\dot{H}^0)}
    \leq  C(h^\gamma+k^{\frac{\gamma}{2}})(\|\varphi\|_{L^2(\Omega;\dot{H}^{\gamma} ) }
         +
         \|\psi\|_{L^2(\Omega;\dot{H}^{\gamma})}).
\end{align}
Before treating $ P_2 \mathbb{J}_2 $, we again recall the stability property \eqref{DDPf5}.
Following the same arguments as used in (\ref{aldlalg}) and using Assumption \ref{ass:F},
(\ref{eq:u(t)-u(s).final}) with $\varrho=\gamma-1$ and (\ref{eq:error-u-U-general}) with $\beta=\gamma$, $i=0$ give
\begin{align}\label{dingliad2}
\| P_2 \mathbb{J}_2 \|_{ L^2(\Omega;\dot{H}^0) }
\leq &
\sum_{j=0}^{n-1} \int_{t_j}^{t_{j+1}} \| P_2 r(k\mathcal{A}_h)^{n-j} \mathbf{P}_h(\mathbf{F}(X^{j})
      -\mathbf{F}(X(s)))\|_{L^2(\Omega;\dot{H}^0)} \, \text{d} s
\nonumber\\
\leq &
    C\sum_{j=0}^{n-1} \int_{t_j}^{t_{j+1}}  \|F(U^j)-F(u(s))\|_{L^2(\Omega;\dot{H}^0)} \, \text{d} s
\nonumber\\
\leq&
C\sum_{j=0}^{n-1}
   \int_{t_j}^{t_{j+1}}  \left(   \|U^j-u(t_j)\|_{L^2(\Omega;\dot{H}^0)}
   +
      \|u(t_j)-u(s)\|_{L^2(\Omega;\dot{H}^0)}    \right) \,\dd s
      \nonumber\\
\leq&
C(h^{\gamma} + k^{ \frac{\gamma}{2} } )
 \big(   1+\|\varphi\|_{L^2(\Omega;\dot{H}^{\gamma})} + \|\psi\|_{L^2(\Omega; \dot{H}^{\gamma-1})} \big).
\end{align}
Similarly as in (\ref{leq31}), we utilize (\ref{eq:error-full-vv}) with $s=0, q=1$, \eqref{eq:F-moment-bound},
and  Assumption \ref{ass:F} to achieve
\begin{align}\label{dingliad3}
\| P_2 \mathbb{J}_3 \|_{L^2(\Omega;\dot{H}^0)}
&\leq
\sum_{j=0}^{n-1} \int_{t_{j}}^{t_{j+1}}
\|P_2 F_{kh}^{n-j}  \mathbf{F}(X(s))\|_{L^2(\Omega;\dot{H}^0)} \, \dd s
\nonumber \\
& \leq C ( h + k^{\frac12} ) \sum_{j=1}^{n-1} k t_{n-j}^{-\frac{1}{2}} \sup_{ s \in[0, T] } \| F(u(s)) \|_{L^2(\Omega;\dot{H}^0)}
\nonumber \\
& \leq C (h + k^{\frac12}) \big( 1+\|\varphi\|_{L^2(\Omega;\dot{H}^{0})} + \|\psi\|_{L^2(\Omega; \dot{H}^{-2})} \big).
\end{align}
To handle the term $ P_2 \mathbb{J}_4 $, by (\ref{eq:linear-deter-v-time-regularity-S}) with $\nu=\frac{1}{2}$  and \eqref{eq:F-moment-bound} one can deduce
\begin{align}\label{dingliad4}
\| P_2 \mathbb{J}_4 \|_{L^2(\Omega;\dot{H}^0)}
 & \leq
  \sum_{j=0}^{n-1}   \int_{t_{j}}^{t_{j+1}}
  \|  P_2(\mathcal{S}(t_{n-j})-\mathcal{S}(t_n-s)) \mathbf{F}(X(s)) \|_{L^2(\Omega;\dot{H}^0)}  \,\dd s\nonumber\\
  &\leq C  \int_0^{t_n} k^{\frac{1}{2}}  (t_n-s)^{-\frac{1}{2}}\,\dd s
  \sup_{ s \in[0, T] } \| F(u(s)) \|_{L^2(\Omega;\dot{H}^0)}\nonumber\\
  &\leq Ck^{\frac{1}{2}} \big( 1+\|\varphi\|_{L^2(\Omega;\dot{H}^{0})} + \|\psi\|_{L^2(\Omega; \dot{H}^{-2})} \big).
\end{align}
Finally, we use  It\^{o}'s isometry, Lemma \ref{lemmad} with $\mu = \gamma$ and  (\ref{eq:error-full-u-v})
with $\beta=\gamma$ to show
\begin{align}\label{dingliad6}
\|  P_2 \mathbb{J}_5 \|^2_{L^2(\Omega;\dot{H}^0)}
= &
\sum_{j=0}^{n-1}\int_{t_{j}}^{t_{j+1}}  \left\|   P_2\left(r(k\mathcal{A}_h)^{n-j}\mathbf{P}_h
-
  \mathcal{S}(t_n-s)\right) \mathbf{B}Q^{\frac{1}{2}} \right\|^2_{\mathrm{HS}}  \dd s
\nonumber \\
\leq &
2 \sum_{j=0}^{n-1} \int_{t_j}^{t_{j+1}} \left(\|P_2 F_{kh}^{n-j} \mathbf{B}Q^{\frac{1}{2}}\|_{\mathrm{HS}}^2
+
\|P_2(\mathcal{S}(t_{n-j}) - \mathcal{S}(t_n-s)) \mathbf{B}Q^{\frac{1}{2}}\|_{\mathrm{HS}}^2 \right) \dd s
\nonumber\\
\leq &
C(h^{2\gamma}+k^{\gamma}) \|A^{\frac{\gamma-1}{2}}Q^{\frac{1}{2}}\|^2_{\mathrm{HS}}.
\end{align}
Gathering (\ref{l12es})-(\ref{dingliad6}) together implies (\ref{qudddansq1}) and the proof of this theorem is thus complete.
$\square$

\section{Numerical examples}
\label{sec:numerical-exam}

\begin{figure}[htp]
      \centering
      \includegraphics[width=7in,height=2in] {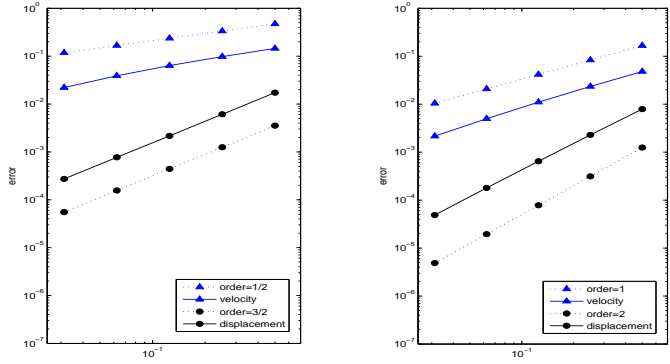}
         \caption{ Mean-square convergence rates for the spatial discretizations (Left: $Q = I$; right: $Q = A^{-0.5005}$)}
         \label{fig:spatial.error}
\end{figure}

In this section, we report some numerical experiments to illustrate our previous findings.
Let us consider the following strongly damped wave equation,
subject to a perturbation of additive noise,
\begin{equation}\label{eq:numer-exam1}
\left\{
    \begin{array}{lll}
    u_{tt} = \Delta u + \Delta u_{t} -\sin( u ) + \dot{W}(t), \quad  t \in (0, 1], \:\: x \in (0,1),
    \\
     u(0, x) = \frac{\partial u}{\partial t}  (0,x) = 0, \: x \in (0,1),
     \\
     u(t, 0) = u(t,1) = 0,  \: t >0.
    \end{array}\right.
\end{equation}
In the following experiments, we aim to test mean-square approximation errors as theoretically measured
in \eqref{thuh1}, \eqref{thuh2}, \eqref{quansq1} and  \eqref{qudddansq1}.  The expectation is
approximated by the Monte-Carlo approximation, using $M = 100$ path simulations. As the first task,
we examine the spatial approximation errors  $\|u(T)-u_h(T)\|_{L^2(\Omega;\dot{H}^0)} $ and
$\|u_t(T)-u_{h,t}(T)\|_{L^2(\Omega;\dot{H}^0)}$, with the endpoint $T = 1$ fixed.
The ``true'' solutions $u(T)$, $u_t(T)$ are identified  with numerical ones using small step-sizes
$ h_{\text{exact}} = 2^{-8}, k_{\text{exact}} = 2^{-14}$. The numerical approximations under various spatial mesh sizes $h=2^{-i}, i =1,2,..., 5$ are achieved via time-stepping with $k_{\text{exact}} = 2^{-14}$.
The resulting computational errors are listed in Table \ref{table:spatial-error}, where two kinds of noises
are considered  including the space-time white noise case ($Q = I$) and the trace-class noise case
($ Q = A^{-0.5005} $). To clearly see the convergence rates, we depict in Figure \ref{fig:spatial.error}
the  errors versus mesh sizes in logarithmic scale. As expected, the slopes of the errors (solid lines) and those of the reference dashed lines match well. More formally, the finite element spatial approximation errors in the space-time white noise case ($Q = I$) exhibit convergence rates of order $\tfrac32$ for the displacement and order $\tfrac12$ for the velocity, which coincides with our previous theoretical findings in Theorem \ref{thm:conteh1} when $\gamma < \tfrac12$.  For the other case ($ Q = A^{-0.5005} $ and $\gamma =1$), the errors show the predicted rates of order $2$  for the displacement and order $1$ for the velocity (see the right plot in Figure \ref{fig:spatial.error}).

\begin{table}
 \centering
\caption{Mean-square spatial errors for the displacement  and the velocity}
\vspace{5mm}
\begin{tabular}{|c|c|c|c|}\hline
 &$ h$ & $ (\mathbf{E}\|u(T) - u_h (T) \|^2)^{\frac12}$ & $(\mathbf{E}\| u_t(T) - u_{h, t}(T) \|^2)^{\frac12}$
 \\\hline
         &$1/2$   &0.017262              & 0.144681       \\\cline{2-4}

        &$1/4$    &0.006098              &0.097764        \\\cline{2-4}

 $Q = I$

         &$1/8$    &0.002158              &0.063434       \\\cline{2-4}

         &$1/16$   &7.694527e-004         &0.038866       \\\cline{2-4}

         &$1/32$   &2.723050e-004         &0.022045         \\\cline{2-4}
\hline

          &$1/2$   &0.007918              & 0.048106             \\\cline{2-4}

         &$1/4$    &0.002289               &0.023401         \\\cline{2-4}

$Q = A^{-0.5005}$

         & $1/8$   & 6.467743e-004         &0.011036          \\\cline{2-4}

        &$1/16$    & 1.800250e-004         &0.004982           \\\cline{2-4}
        &$1/32$     &4.888214e-005         &0.002160           \\\cline{2-4}
\hline
\end{tabular}
\label{table:spatial-error}
\end{table}

\begin{figure}[htp]
    \centering
    \includegraphics[width=7in,height=2in] {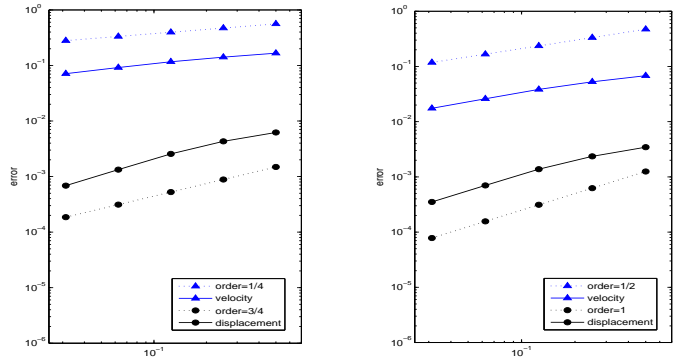}
         \caption{ Mean-square convergence rates for the temporal discretizations (Left: $Q = I$; right: $Q = A^{-0.5005}$)}
         \label{fig:time.error}
\end{figure}

Next, we proceed to tests on the convergence rates of temporal approximations. To this end, we fix
$h = 2^{-7}, T = 1 $ and measure $ \| u_h (T) - U^N \|_{L^2(\Omega;\dot{H}^0)}$ and $ \| u_{h,t} (T) - V^N \|_{L^2(\Omega;\dot{H}^0)}$ for five different time stepsizes $k = \tfrac{1}{N}, N = 2^{3}, 2^{4},..., 2^{7} $.
In order to obtain $u_h (T), u_{h,t} (T)$,  we perform time-stepping using small time stepsize
$k_{\text{exact}} = 2^{-12}$. In Table \ref{table:temporal-error} we present the computational errors for the two noise cases $Q = I$
and $Q = A^{-0.5005}$.  Similarly as before, these approximation errors are plotted versus time step-sizes in Figure \ref{fig:time.error},
where one can easily observe the expected convergence rates. For example,  in the trace-class noise case when
$Q = A^{-0.5005}$, the approximation errors for the displacement  and the velocity decrease at slopes of $1$ and $\tfrac12$, respectively, as the time step-sizes decrease.  Also, convergence rates of order $\tfrac34$ and $\tfrac14$ are detected for the displacement  and the velocity in the space-time white noise case (see the left plot in Figure \ref{fig:time.error}). All in all, the above observations are all consistent with the previous theoretical results.

\begin{table}
 \centering
\caption{Mean-square temporal errors for the displacement  and the velocity}
\vspace{5mm}
\begin{tabular}{|c|c|c|c|c|c|}\hline
 &$k$&$(\mathbf{E}\|u_h (T) - U^N\|^2)^{\frac12}$&$(\mathbf{E}\| u_{h, t} (T) - V^N \|^2)^{\frac12}$
 \\\hline
         &$1/8$    &0.006226              & 0.166427           \\\cline{2-4}

         &$1/16$   &0.004302             &0.141315      \\\cline{2-4}

 $Q = I$

         &$1/32$   &0.002560               &0.116514        \\\cline{2-4}

         &$1/64$   &0.001332               &0.091829        \\\cline{2-4}

         &$1/128$  &6.853130e-004          &0.071157       \\\cline{2-4}
\hline

         &$1/8$    &0.003446               &0.068094       \\\cline{2-4}

        &$1/16$     & 0.002356             &0.052651       \\\cline{2-4}

$Q = A^{-0.5005}$

        &$1/32$     &0.001377               &0.038405        \\\cline{2-4}
        &$1/64$     &6.993377e-004          &0.025994        \\\cline{2-4}
        &$1/128$    &3.512776e-004          &0.017476         \\\cline{2-4}
\hline
\end{tabular}
\label{table:temporal-error}
\end{table}

\section*{Acknowledgments}
Part of this work was done when the authors attended a seminar on stochastic computations at AMSS, Beijing.
The authors want to thank Prof. Jialin Hong for his kindness and help during their stay.
\bibliographystyle{abbrv}

\bibliography{bibfile}

 \end{document}